\documentclass[12pt,a4paper]{amsart}

\usepackage{amssymb}

\makeatletter
\renewcommand{\@begintheorem}[2]{
\rm \trivlist \item [\hskip \labelsep {\bf #2\ \ #1.}]
                                }
\makeatother

\makeatletter

\DeclareFontFamily{U}{cyr}{}
\DeclareFontShape{U}{cyr}{m}{n}{
  <5> wncyr5 <6> wncyr6 <7> wncyr7 <8> wncyr8 <9> wncyr9 <10->
wncyr10}{}
\DeclareMathAlphabet{\mathcyr}{U}{cyr}{m}{n}

\input cyracc.def

\newcommand{\ts}{\vspace{\baselineskip}\noindent{\bf Proof.}$\;\;$}
\newcommand{\ZZ}{{\bf Z}}
\newcommand{\QQ}{{\bf Q}}
\newcommand{\RR}{{\bf R}}
\newcommand{\CC}{{\bf C}}

\newcommand{\FF}{{\bf F}}

\newcommand{\PP}{{\bf P}}

\newcommand{\bA}{{\mathbb A}}

\newcommand{\bP}{{\mathbb P}}

\newcommand{\cL}{{\mathbb L}}

\newcommand{\cO}{{\mathcal O}}

\newcommand{\hY}{{\widetilde{Y}}}
\newcommand{\hZ}{{\widetilde{Z}}}

\newcommand{\et}{{\text{\'et}}}

\newcommand{\rmd}{\mbox{d}}
\newcommand{\gop}{{\wp}}
\newcommand{\bes}{\begin{equation*}}
\newcommand{\ees}{\end{equation*}}

\title{Geometry and arithmetic of Maschke's Calabi-Yau threefold}
\author{Gilberto Bini}
\author{Bert van Geemen}
\address{Dipartimento di Matematica, Universit\`a di Milano,
Via Saldini 50, I-20133 Milano, Italia}
\email{gilberto.bini@unimi.it}
\email{lambertus.vangeemen@unimi.it}

\begin{document}

\begin{abstract}
Maschke's Calabi-Yau threefold is the double cover of projective three space
branched along Maschke's octic surface. This surface is defined by the lowest
degree invariant of a certain finite group acting on a four dimensional vector space.
Using this group, we show that the middle Betti cohomology group of the threefold decomposes into the direct sum of $150$ two-dimensional Hodge substructures.
We exhibit one dimensional families of rational curves on the threefold and verify that the associated Abel-Jacobi map is non-trivial.
By counting the number of points over finite fields,
we determine the rank of the N\'eron-Severi group of Maschke's surface and
the Galois representation on the transcendental lattice of some of its quotients.
We also formulate precise conjectures on the modularity of the Galois representations
associated to Maschke's threefold and to a genus $33$ curve which parametrizes rational curves in the threefold.
\end{abstract}

\maketitle

The Hodge structure on the middle dimensional Betti cohomology group of a
Calabi-Yau threefold carries important information on the moduli and
the one-dimensional algebraic cycles on the threefold.
However, if the threefold is easy to define, say by one equation in a (weighted) projective space, the dimension $h^3$ of this vector space tends to be large. For example, a smooth quintic threefold in $\PP^4$ has $h^3=204$ and a double octic, i.e.\ a double cover of $\PP^3$ branched along a smooth surface of degree $8$, has $h^3=300$. Using automorphisms of the threefolds, one can decompose the cohomology into subrepresentations, which give rise to Hodge substructures.
In this paper we consider a double octic with a particularly large automorphism group $\overline{G}$ of order $16\cdot (6!)=11520$.

In 1887 Heinrich Maschke studied the invariants of a finite group $G$,
related to genus two theta functions,
acting on a four dimensional complex vector space (\cite{M}).
The group $\overline{G}$ is the image of $G$ in $Aut(\PP^3)$,
it is the quotient of $G$ by its center, which is cyclic of order four.
Maschke's main result is the determination of generators of the ring of $G$-invariants.

In particular, the lowest degree invariant has degree $8$ and is given
by (\cite{M}, (12))
$$
F\,:=\,\sum_{i=0}^3 x_i^8\,+\,14\sum_{i<j}x_i^4x_j^4\,+\,168x_0^2x_1^2x_2^2x_3^2.
$$
The octic surface $S$ in $\PP^3$ defined by this polynomial
$$
S\,:=\,Zeroes(F)\qquad(\subset\PP^3)
$$
is a smooth surface, which we call Maschke's surface.
This surface, and a subgroup $G_8$ of its automorphism group,
appeared recently in  \cite{BS}.
The double cover $X$ of $\PP^3$, branched along $S$, is a smooth Calabi-Yau threefold which we call Maschke's double octic.

Using an easy generalization of a formula of Ch\^enevert \cite{Ch}
(see section \ref{trfo})
and a computer, we find in section \ref{decXG}
that $H^3(X,\QQ)$ is a direct sum of two dimensional Hodge substructures.
As the action of $G$ on $\CC^4$ is given by matrices with coefficients in the field $\QQ(i)$, $i^2=-1$, the \'etale cohomology group $H^3_\et(X,\QQ_\ell)$ is decomposed into two dimensional
$G_{\QQ(i)}:=Gal(\overline{\QQ}/\QQ(i))$-representations.
Our numerical results,
using point counting over finite fields with MAGMA (\cite{magma})
and the Lefschetz fixed point formula,
suggest that $H^3_\et(X,\QQ_\ell)$ is a direct sum of two dimensional $G_{\QQ}:=Gal(\overline{\QQ}/\QQ)$-representations.
For such representations there are now
various modularity results, and we did succeed in matching the numerical
data to explicitly given elliptic modular forms.
However, we were not able to provide a proof for our conjectured decomposition
given in section \ref{galX}.

The decomposition of $H^3(X,\QQ)$ shows that the
Griffiths intermediate Jacobian of $X$ has an abelian subvariety $J(X)_a$
of dimension $149$.
In section \ref{lines} we describe a family of rational curves on $X$,
parametrized by a curve $\widetilde{C}_+$, such that the Abel-Jacobi map
$J(\widetilde{C}_+)\rightarrow J(X)_a$ is non-constant. This provides some evidence for
the generalized Hodge conjecture for $X$. We also analyze the \'etale cohomology of $\widetilde{C}_+$ in section \ref{arcp}. It seems that $J(\widetilde{C}_+)$ is isogeneous
to a product of $33$ elliptic curves, and that the image of the Abel-Jacobi map
might be `as big as possible', that is, that the image might have dimension $24$.

A somewhat similar study of the Abel-Jacobi map for complete intersections
of four quadrics in $\PP^7$ with an involution was done in \cite{Fabio}.
The action of a finite group and the Galois representations associated with
the quintic threefolds of the Dwork pencil were analyzed in \cite{COR}.
In that case, each threefold has a one dimensional family of lines, which was analyzed in \cite{Anca}, and recently Candelas obtained an explicit description of this family \cite{Cp}.
As far as we know, the Abel-Jacobi map has not yet been described in detail.
The paper \cite{HV} considers surfaces which are products of a $\PP^1$ and
an elliptic curve $E$ in certain Calabi-Yau threefolds and
shows that in some cases $H^1(E,\QQ)\cong H^3(E\times \PP^1,\QQ)$ contributes to $H^3(X,\QQ)$.
This can be viewed as an instance of an injective Abel-Jacobi map associated to a family of $\PP^1$'s parametrized by $E$.

The group $G$ also acts on Maschke's surface $S$ and in sections
\ref{cohS}, \ref{quotS} and \ref{H2etS}
we use it to study the Hodge structure on $H^2(S,\QQ)$ and
the Galois representation on $H^2_\et(S,\QQ_\ell)$.
This allows us to determine the rank of the N\'eron-Severi
group of $S$,
from which we conclude that the classes of lines in $S$ generate
a subgroup of finite index.
We also give projective models of some quotients of $S$.

In the first section
we briefly recall how a finite group of automorphisms of an
algebraic variety $X$
can be used to decompose the cohomology of $X$.
The group $G$ is introduced in section \ref{grG}.

We are indebted to Matthias Sch\"utt for his comments on a first draft of this
paper.

\section{Decomposing cohomology groups}\label{dcg}

\subsection{} The Betti cohomology groups $H^k(X,\QQ)$ of a smooth complex projective variety $X$
are rational Hodge structures of weight $k$.
A subgroup $G$ of $Aut(X)$ induces a decomposition of these groups into Hodge substructures. Similarly, if $X$ and $G$ are defined over $\QQ$, one obtains a
corresponding (via the comparison theorem) decomposition of the \'etale cohomology groups of $X$ into $G_\QQ$-representations.

We briefly recall the basics of rational Hodge structures
(see for example \cite{V1}, Chapter 7), representations of finite groups and their applications to algebraic geometry.

\subsection{Rational Hodge structures}\label{rhs}
The Betti cohomology groups $H^k(X,\QQ)$ of a smooth complex projective variety $X$
are rational Hodge structures of weight $k$.
A rational Hodge structure is a finite dimensional $\QQ$-vector space with a decomposition
of its complexification
$V_\CC:=V\otimes\CC=\oplus V^{p,q}$ with $p,q\in\ZZ_{\geq 0}$, $p+q=k$ and
$\overline{V^{p,q}}=V^{q,p}$.
A morphism of Hodge structures is a $\QQ$-linear map $f:V\rightarrow W$ such that its $\CC$-linear extension satisfies $f(V^{p,q})\subset W^{p,q}$ for all $p,q$. Morphisms between varieties induce morphisms of Hodge structures on the Betti cohomology groups. For a Hodge structure
$V$, the (endo)morphisms of Hodge structures $V\rightarrow V$ form a $\QQ$-algebra (so $\lambda f+\mu g$, $fg:=f\circ g$ are morphisms of Hodge structures if $f,g$ are, for any $\lambda,\mu\in\QQ$).
Thus if $G\subset Aut(X)$ is a subgroup, we get an action, by Hodge endomorphisms, of $G$ on each $H^k(X,\QQ)$.

A subspace $W\hookrightarrow V$ is called a Hodge substructure if
$W_\CC=\oplus W^{p,q}$ where $W^{p,q}:=W\cap V^{p,q}$, in that case $W$ is a Hodge
structure and the inclusion is a morphism of Hodge structures.
The kernel and image of a morphism of Hodge structures are easily seen to be Hodge substructures.

A Hodge structure $V$ is called simple if the only Hodge substructures of $V$ are $\{0\}$ and $V$.
As the Hodge structure on $H^k(X,\QQ)$ is polarized, it is a direct sum of simple
Hodge structures and these are unique up to isomorphism.

An automorphism $\phi:X\rightarrow X$ induces endomorphisms of Hodge structures $\phi^*$ on $H^k(X,\QQ)$ for each $k$. For any $a_0,\ldots,a_m\in \QQ$, one has the Hodge endomorphism $\sum_{i=0}^m a_i(\phi^*)^i$,
the kernel and image of which are thus Hodge substructures of $V$.

\subsection{Representations of finite groups}
We recall some of the basics on representations of finite groups.
Let $G$ be a finite group, let $G=C_1\coprod \ldots\coprod C_M$ be the partition of
$G$ in conjugacy
classes (so $g,g'\in C_i$ for some $i$ iff $g=hg'h^{-1}$ for some $h\in G$).
A representation of $G$ is a homomorphism of groups $\rho:G\rightarrow GL(V)$
where $V$ is a finite dimensional vector space over a field $K$.
A representation is said to be irreducible (over $K$) if the only subspaces $W\subset V$ with $\rho(g)W\subset W$ for all $g\in G$ are $W=\{0\}$ and $W=V$.
A representation $\rho':G\rightarrow GL(V')$ is said to be equivalent (or isomorphic)
to $\rho$ if there is a $K$-linear isomorphism $A:V\rightarrow V'$ such that
$\rho'(g)=A\rho(g)A^{-1}$ for all $g\in G$.

The group $G$ has $M$ irreducible complex (i.e.\ $K=\CC$) representations,
we denote them by $\rho_1,\ldots,\rho_M$ on vector spaces $V_1,\ldots,V_M$.
A complex representation of $G$ (so $K=\CC$) is isomorphic to a direct sum of the irreducible representations, $V\cong \oplus V_i^{n_i}$.

\subsection{Characters}
The isomorphism class of a complex representation $\rho$ is determined by
its character $\chi_\rho$, so $\rho,\rho'$ are isomorphic representations iff $\chi_\rho(g)=\chi_{\rho'}(g)$ for all $g\in G$.
Here the character of $\rho$ is the function
$$
\chi_\rho:\,G\,\longrightarrow\, \CC, \qquad \chi_\rho(g)\,:=\,Tr(\rho(g)),
$$
where $Tr$ is the trace of a linear map. This function is constant on the conjugacy classes. Note that if $V\cong \oplus V_i^{n_i}$, then $\chi_\rho=\sum n_i\chi_i$ where
$\chi_i:=\chi_{\rho_i}$ is the character of the irreducible representation $\rho_i$.
The decomposition of a character into irreducibles is easily done using a scalar product on the space of characters for which the $\chi_i$ are an orthonormal basis.
The trace of the action of an automorphism
on the primitive cohomology of a hypersurface or a ramified cover of $\PP^{n+1}$
can be found in Proposition \ref{propch}.

 \subsection{Decomposing representations}\label{decrep}
Given a complex representation $\rho$ of $G$ and a conjugacy class $C_j$,
consider the linear map
$$
c_{j,\rho}\,:=\,\sum_{g\in C_i}\rho(g)\qquad(\in End(V)).
$$
By construction, $\rho(h)c_{j,\rho}\rho(h)^{-1}=c_{j,\rho}$ for all $h\in G$.
By Schur's lemma, $c_{j,\rho}$ is then scalar multiplication by a complex number $\lambda_{i,j}$ on any irreducible subrepresentation $V_i\subset V$.
This scalar can be easily found by computing the trace of $c_{j,\rho}$ on $V_i$
in two ways. The first is simply $Tr(c_{j,\rho}|V_i)=(\dim V_i)\lambda_{i,j}$, the other
uses a character:
$
Tr(c_{j,\rho}|V_i)=(\sharp C_j)\chi_i(g),
$
for any $g\in C_i$. Thus on the subspace $V_i^{n_i}$ of $V$ we get, for any $g\in C_i$:
$$
\lambda_{i,j}\,:=\,(c_{j,\rho})_{|V_i^{n_i}}\,=\,\frac{\sharp C_j}{\dim V_i}\chi_i(g),\qquad\mbox{define}\quad
\widetilde{V}_{i}:=\cap_{j=1}^M\,\ker(c_{j,\rho}\,-\,\lambda_{i,j}).
$$
Then, by construction, $\widetilde{V}_{i}$ is a subrepresentation of $V$, we have
$V_i^{n_i}\,\subseteq\,\widetilde{V}_i$.
One actually has equality here, since  $Tr(\rho(g)_{|\widetilde{V}_{i}})=n_i\chi_i(g)$, for all $g\in G$ which implies
that the character of the representation $\rho(g)_{|\widetilde{V}_{i}}$
coincides with the character of $V_i^{n_i}$. Hence these representations are isomorphic:
$$
\widetilde{V}_i\,\cong\,V_i^{n_i},
$$
and this provides an explicit method to decompose $V$ into isotypical components.

\subsection{Applications}
In this paper we use the action of finite groups $G\subset Aut(X)$ of automorphisms of smooth complex varieties on the Betti cohomology groups $H^k(X,\QQ)$ of $X$.
In particular, if $\rho_i$ is an irreducible representation of $G$ such that
its character satisfies $\chi_i(g)\in\QQ$ for all $g\in G$, then each $\lambda_{i,j}\in\QQ$ and we can split
$$
H^k(X,\QQ)\,=\,\widetilde{V}_i\oplus \widetilde{V}_i'\qquad
\widetilde{V}_i\,:=\,\cap_{j=1}^M\,\ker(c_{j,\rho}\,-\,\lambda_{i,j}),\qquad(\chi_i(g)\in\QQ\quad\forall g\in G),
$$
where $\rho:G\rightarrow GL(H^k(X,\QQ))$ is the action of $G$,
$\rho(g):=g^*$.
As observed in section \ref{rhs}, $\widetilde{V}_i$ and $\widetilde{V}_i'$ will be rational Hodge structures.

Similarly, one can decompose the \'etale cohomology group
$H^k_{\et}(X,\QQ_\ell)=\widetilde{V}_{i,\ell}\oplus \widetilde{V}'_{i,\ell}$. The comparison isomorphism between \'etale and Betti cohomology groups implies that these decompositions coincide after embedding $\QQ_\ell$ into $\CC$ and tensoring by $\CC$.

\subsection{Commuting group actions}\label{cga}
Let again $\rho:G\rightarrow GL(V)$ be a representation of the finite group $G$ on a $K$-vector space $V$ and let
$$
V\,\cong\,\oplus_i \widetilde{V}_i,\qquad \widetilde{V}_i\,\cong\,V_i^{n_i}
$$
be the decomposition into representations $V_i$ which are irreducible over $K$.
More intrinsically, one has
$$
\widetilde{V}_i\,\cong\,V_i\otimes W_i,\qquad W_i\,:=\,Hom_G(V_i,V)\,=\,Hom_G(V_i,\widetilde{V}_i)
\cong K^{n_i},
$$
the $G$-equivariant maps; the isomorphism is given by $v\otimes f\mapsto f(v)$.
In case another group $H$ acts linearly on $V$ and the actions of $G$ and $H$ commute,
one has an action of the group $G\times H$ on $V$.
Essentially by Schur's lemma (cf.\ \cite{GW}, $\S$3.3) each $W_i$ is then
an $H$-representation space and the action of $G\times H$ on $V$ is given by
$$
V\,=\,\oplus_i V_i\otimes W_i,\qquad (g,h)(v\otimes w)\,=\,(gv)\otimes (hw).
$$
In particular, as an $H$-representation one has
$$
V\,\cong\,\oplus_iW_i^{m_i},\qquad m_i:\,=\,\dim V_i.
$$

In case $G$ acts on $V=H^k(X,\QQ)$ by morphisms of Hodge structures, then $W_i$ is a rational Hodge structure
and
$\widetilde{V}_i\cong W_i^{m_i}$, where $m_i=\dim V_i$, is an isomorphism
of rational Hodge structures.
In fact, the rational Hodge structure on $\widetilde{V}_i$ is defined by a representation $$
h\,:\,\CC^*\,\longrightarrow\,GL(\widetilde{V}_{i,\RR}),
\qquad h(z)v_{p,q}:=z^p\overline{z}^q v_{p,q}
$$
where $v_{p,q}\in (\widetilde{V}_i)^{p,q}$ and the action of $h$ is extended
$\CC$-linearly to $\widetilde{V}_{i,\CC}$ where
$\widetilde{V}_{i,K}:=\widetilde{V}_i\otimes_\QQ K$.
As $g\in G$ preserves the Hodge structure, it must commute with $h(z)$ for all $z$ and hence $G\times \CC^*$ acts on $V_{i,\RR}$.
The representation $h$ is then obtained from a representation
$h':\CC^*\rightarrow GL(W_i)$. This representation defines the Hodge structure
on $W_i$.

In case $X\subset\PP^n$ is defined by equations with coefficients in a number field $L$ and the action of any $g\in G\subset Aut(X)$ is given by maps with coefficients in
$L$, the maps $g^*$ will commute with the action of $G_L:=Gal(\overline{L}/L)$ on
$V_{\ell}:=H^k_{\et}(X,\QQ_\ell)$. Thus, as a $G\times G_L$-representation,
$$
V_{\ell}\,\cong\, \oplus_i \widetilde{V}_{i,\ell}\,\cong\,
\oplus_iV_{i,\ell}\otimes W_{i,\ell}.
$$
In particular, the $G$-representation $\widetilde{V}_{i,\ell}\cong V_{i,\ell}^{n_i}$ defines a Galois representation
on $W_{i,\ell}$, with
$\dim W_{i,\ell}=n_i$, and as $G_L$-representations $\widetilde{V}_{i,\ell}\cong
W_{i,\ell}^{m_i}$ where $m_i:=\dim V_i$.

\subsection{Cycles and Motives}
The splitting of the various cohomology groups comes from cycles on the product $X\times X$ and is thus motivic. The graph $\Gamma_g:=\{(x,gx):x\in X\}$ has a class $[\Gamma_g]\in H^{2n}(X\times X)$ with $n=\dim_\CC X$ (where one can choose any Weil cohomology group).
Using the K\"unneth theorem and Poincar\'e duality to identify $End(H^k(X))$ with $H^{2n-k}(X)\otimes H^{k}(X)$, the K\"unneth components of this class induce maps, compatible with the Hodge structure or the Galois action, which are just the $g^*$ on $H^k(X)$:
$$
g^*\,=\,[\Gamma_g]_k:\,H^k(X)\,\longrightarrow \, H^k(X).
$$
The operators $c_{j,\rho}$ are thus induced from (the action of) certain cycles in $CH^*(X\times X)$ and will act on any Weil cohomology group of $X$.

\section{The group $G$} \label{grG}

\subsection{The Heisenberg group $H$}\label{hgrp}
For the groups in these sections we refer to \cite{CG}.

The Heisenberg group is defined as the set
$$
H=\mu_4 \times (\ZZ/2\ZZ)^2\times (\ZZ/2\ZZ)^2,\qquad
\mu_4:=\{z\in\CC:\,z^4=1\,\},
$$
so, with the obvious group structure, $H$ would be isomorphic to $(\ZZ/4\ZZ)\times (\ZZ/2\ZZ)^4$. However, we define
the group operation by
$$
(s,x,x^*)(t,y,y^*)\,:=\,(st(-1)^{y^*(x)},x+y,x^*+y^*)
\quad  \mbox{with}\quad y^*(x):=y_1^*x_1+y_2^*x_2,
$$
where $x=(x_1,x_2),y^*=(y_1^*,y_2^*)\in(\ZZ/2\ZZ)^2$.
One easily verifies that
$$
(s,x,x^*)^{-1}=(s^{-1}(-1)^{x^*(x)},x,x^*)
$$
and that the commutator of two elements in $H$ is given by:
$$
(s,x,x^*)(t,y,y^*)(s,x,x^*)^{-1}(t,y,y^*)^{-1}=((-1)^{x^*(y)-y^*(x)},0,0).
$$
The center of $H$ is $\mu_4$ and the commutator defines a symplectic form $E$ on
the $\FF_2$-vector space $H/\mu_4\cong (\ZZ/2\ZZ)^4$:
$$
E((x,x^*),(y,y^*))\,=\,y^*(x)-x^*(y)=y_1^*x_1+y_2^*x_2+x_1^*y_1+x_2^*y_2,
$$
with $x=(x_1,x_2)$, $x^*=(x_1^*,x_2^*)$, $y=(y_1,y_2)$ and
$y^*=(y_1^*,y_2^*)$, all in $(\ZZ/2\ZZ)^2$.

\subsection{The Schr\"odinger representation}
There is a (unique, faithful) representation $U$ of the finite group $H$ on the vector space
$\CC^4$, called the Schr\"odinger representation, such that $s\in\mu_4$ acts as scalar multiplication by $s$.
Identifying $\CC^4$ with the $\CC$-vector space of maps $f:(\ZZ/2\ZZ)^2\rightarrow\CC$
this representation is defined as:
$$
\Bigl(U_{(s,x,x^*)} f\Bigr)(z) \,:=\,s(-1)^{x^*(z)} f(x+z),
$$
where $U_{(s,x,x^*)}\in GL(4,\CC)$ gives the action of $(s,x,x^*)\in H$.
The $\delta$-functions provide a basis of this vector space: $\delta_x(y)=0$ if $x\neq y$ and $\delta_x(x)=1$ for $x,y\in (\ZZ/2\ZZ)^2$.
The Schr\"odinger representation is then:
$$
U_{(s,x,x^*)}\delta_a\,:=\,s(-1)^{x^*(x+a)}\delta_{x+a}.
$$

Identifying $\{0,1\}$ with $\ZZ/2\ZZ$ and writing
$$
x_0=\delta_{00},\quad x_1=\delta_{01},\quad x_2=\delta_{10},\quad x_3=\delta_{11}
$$
for the basis of $\CC^4$ and writing
$$
h_{abcd}\,:=\,U_{(1,(a,b),(c,d))},\qquad c\,:=\,U_{(i,(0,0),(0,0))},\quad i^2\,=\,-1,
$$
the
Schr\"odinger representation on the 4 variables $x_0,\ldots,x_3$ is:
$$
x:=(x_0,x_1,x_2,x_3)\,\longmapsto\left\{
\begin{array}{rcl}
h_{0001}(x)&=&(x_0,-x_1,x_2,-x_3)\\
h_{0010}(x)&=&(x_0,x_1,-x_2,-x_3)\\
h_{0100}(x)&=&(x_1,x_0,x_3,x_2)\\
h_{1000}(x)&=&(x_2,x_3,x_0,x_1)\\
c(x)&=&(ix_0,ix_1,ix_2,ix_3).
\end{array}\right.
$$

\subsection{The group $G$}\label{mgrp}
For convenience, we will now identify the abstract group $H$
with its image under $U$ in $GL(4,\CC)$.
The normalizer $N$ in $GL(4,\CC)$ of $H$ is defined as:
$$
N\,:=\,\{\,M\in\, GL(4,\CC)\,:\,MHM^{-1}=H\,\}.
$$
An element $M\in N$ induces an automorphism of $H$
which we also denote by $M$,
$$
M:\,H\,\longrightarrow\, H,\qquad
h\,\longmapsto\,{h'}\qquad\mbox{if}\quad
MU_hM^{-1}\,=\,U_{h'}.
$$

As $N$ acts by automorphisms on $H$, each $M\in N$ induces a linear map
$\phi_M$ on the quotient of $H$ by its center, $H/\mu_4\cong (\ZZ/2\ZZ)^4$.
These linear maps are easily seen to preserve the symplectic form $E$:
$$
E(\phi_M v,\phi_M w)\,=\,E(v,w),\qquad
(M\in N,\;\phi_M\in Aut((\ZZ/2\ZZ)^4)=GL(4,\FF_2)),
$$
for all $v,w\in (\ZZ/2\ZZ)^4$.
Thus we get a homomorphism to a finite symplectic group:
$$
N\,\longrightarrow\,Sp(4,\FF_2),\qquad M\,\longmapsto\, \phi_M
$$
which can be shown to be surjective.
The group $Sp(4,\FF_2)$ is isomorphic to the symmetric group $S_6$: see, for instance, \cite{CG}, Appendix C.

Any element $U_w\in H$ is obviously an element of $N$.
Since $U_wU_vU_w^{-1}=(-1)^{E(v,w)}U_v$
we get $\phi_M=I$ if $M=U_w\in H\subset N$.
The homomorphism above fits in an exact sequence:
$$
0\,\longrightarrow\,\CC^\times\cdot H\,\longrightarrow\,N\,\longrightarrow\,
S_6\,\longrightarrow 0.
$$

The group $G$ is a subgroup of $N$ which maps onto $S_6$ and such that
$G\cap (\CC^\times\cdot H)\,=\,H$. It thus has order
$$
\sharp G\,=\,(\sharp H)(\sharp S_6)\,=\,2^6\cdot (6!)\,=\,46080.
$$
The group $G$ is generated by the following two elements in $GL(4,\CC)$:
$$
g_1\,:=\,\begin{pmatrix} 1&0&0&0\\0&1&0&0\\0&0&i&0\\0&0&0&i\end{pmatrix},
\qquad
g_2:=\frac{1}{2}
\begin{pmatrix} -1&-i&-i&-1\\
i&1&-1&-i\\
i&-1&1&-i\\
1&-i&-i&1
\end{pmatrix},\quad i:=\sqrt{-1}.
$$
The homomorphism $G\rightarrow S_6$ can be chosen so that the images of $g_1,g_2$ in $S_6$ are the cycles $(12)$ and $(123456)$ respectively. We used the computer program MAGMA (\cite{magma}) for computations involving $H$ and $G$.

\section{The second cohomology group of $S$}\label{cohS}

\subsection{}We study the rational Hodge structure and the $G$-action
on $H^2(S,\QQ)$.

\subsection{The Hodge decomposition of $H^2(S,\QQ)$}
The second (singular) cohomology group with coefficients in $\QQ$
of an octic surface in $\PP^3$ has dimension $b_2(S)=302$. As for any smooth  surface,  this rational Hodge structure splits as
$$
H^2(S,\QQ)\,=\,T_{S,\QQ}\,\oplus\,NS(S)_\QQ\qquad
(
NS(S)_\QQ:=NS(S)\otimes_\ZZ\QQ).
$$
The Hodge substructure $NS(S)_\QQ$ is the N\'eron-Severi group of $S$ tensored with $\QQ$, it is pure of type $(1,1)$.
The transcendental substructure,
which is the orthogonal complement of $NS(S)_\QQ$ w.r.t.\ the intersection form on $H^2$, has the Hodge decomposition
$$
T_{S,\QQ}\otimes_\QQ\CC\,=\,T_S^{2,0}\,\oplus\,T_S^{1,1}\,
\oplus\,T_S^{0,2},\qquad
\overline{T_S^{2,0}}\,=\,T_S^{0,2}.
$$
One has $h^{2,0}(S)=\dim T_S^{2,0}=\dim h^0(\omega_S)=p_g(S)$, the geometric genus of $S$ (and $\omega_S$ is the canonical bundle of $S$).
For an octic surface we have, by adjunction, that
$\omega_S\cong \cO_S(4)$ and thus
$$
H^{2,0}(S)\cong R_4,\qquad\mbox{where}\quad R_4=H^0(\PP^3,\cO(4))
$$
is the
complex vector space of homogeneous polynomials of degree four in $x_0,\ldots,x_3$.
Hence $h^{2,0}(S)=\dim R_4=35$ and
$$
b_2(S)\,=\,\dim H^2(S,\QQ)=302,\qquad
\left\{\begin{array}{rccl}
h^{2,0}(S)&=&\dim H^{2,0}(S)&=35,\\
h^{1,1}(S)&=&\dim H^{1,1}(S)&=302-2\cdot 35=232.
\end{array}\right.
$$
From this we get that $\dim NS(S)_\QQ\leq 232$ and $\dim T_{S,\QQ}\geq 70$.
We need to work a little harder to actually find these dimensions, they are
$\dim NS(S)_\QQ=202$, $\dim T_{S,\QQ}=100$, see sections \ref{dimdh2},  \ref{decS}.

\subsection{The N\'eron-Severi  group of $S$}\label{NS(S)}
To get a lower bound on the rank of the N\'eron-Severi group of $S$,
we study the subgroup $L_S$ spanned by the lines in $S$.
In \cite{BS} it was shown that the surface $S$ contains
exactly $352$ lines.
Under the action of $G$, there are two orbits on the set of lines.
The orbit of the line
$$
l_3:\quad \langle \,(\alpha:1:0,0),\;(0:0:\alpha:1)\,\rangle,
\qquad \alpha^4-2\alpha^3+2\alpha^2+2\alpha+1=0,
$$
($\alpha=\zeta_{12}+\zeta_{12}^2$, where $\zeta_{12}$
is a primitive $12$-th root of unity)
contains 160 lines, whereas the orbit of $l_5$ contains 192 elements
with
$$
l_5:\quad\langle(1:a:a\sqrt{-1}:0), \, (0:a:-a\sqrt{-1}:1)  \rangle,\qquad
a:=(1+\sqrt{-1})(1+\sqrt{5})/4.
$$
Using a computer, we found that the rank of the symmetric matrix of intersection numbers
$(l,m)$, where $l,m$ run over the $256$ lines on $S$, is $202$. Hence these lines
span a subspace $L_S$ of dimension $202$ of
$NS(S)_\QQ$ $(\subset H^{1,1}(S))$.
We will see later that actually $L_S=NS(S)_\QQ$, but at this point we only know:
$$
\dim L_S=202\leq \dim NS(S)_\QQ\leq h^{1,1}(S)=232.
$$

\subsection{The action of $G$ on $H^2(S,\QQ)$}\label{dimdh2}
Using Magma, we found the character table of $G$. There are $59$ irreducible representations. Using a formula of Ch\^{e}nevert,
see Proposition \ref{propch},
it is easy to compute the traces of elements of $G$ on $H^2(S,\QQ)_{pr}$, the primitive cohomology, and thus to find the decomposition of the primitive cohomology into irreducible $G$-representations. The decomposition of
$H^2(S,\QQ)$ is then just the sum of this decomposition with the one dimensional trivial representation.

The result is a bit complicated: the primitive cohomology group
$H^2(S,\QQ)_{pr}$ is the direct sum of $10$
distinct irreducible representations,
(which we do not describe explicitly) with:
{\renewcommand{\arraystretch}{1.2}
$$
\begin{array}{|r|r|r|r|r|r|r|r|r|r|r|}\hline
\mbox{dimension}&
1&5& 5& 10& 15& 15& 30& 45& 45& 45\\ \hline
\mbox{multiplicity}&
1&1&3&1&2&1&3&1&1&1\\ \hline
\end{array}
$$
}
Moreover, the characters of these $10$ representations take values in $\ZZ$,
so we can split
$$
H^2(S,\QQ)_{pr}\,=\,\bigoplus_{i=1}^{10}\; V_i\otimes W_i,
$$
where the $V_i$ are irreducible $G$-representations and $W_i$ are rational Hodge structures, and $\dim W_i$ is the multiplicity of $V_i$ in $H^2(S,\QQ)_{pr}$.

To relate this decomposition to the one in the transcendental and trivial
(N\'eron-Severi) Hodge substructures, we used the isomorphism of
$G$-representations
$H^{2,0}(S)\cong R_4(\epsilon)$, where $\epsilon$ is the unique non-trivial one dimensional representation of $G$ (recall that $G/H\cong S_6$, and $\epsilon$ factors over the sign representation of $S_6$).

We found that $H^{2,0}(S)$ is the direct sum of two irreducible representations,
one of dimension 5 and one of dimension 30.
As the $10$ subrepresentations of $H^2(S,\QQ)$ have characters with values in $\ZZ$, the direct sum of these two representations is also isomorphic to $\overline{H^{2,0}}=H^{0,2}$.
In particular, their multiplicity in $H^2(S,\QQ)_{pr}$ is at least $2$.
This suffices to identify them as the $V_i$ with $i=3,7$ in the order given in the table above, so:
$$
H^{2,0}(S)\,\cong\,H^{0,2}(S)\,\cong\,V_3\oplus V_7,\qquad
\dim V_3=5,\quad \dim V_7=30.
$$
Applying the methods of section \ref{cga} to the Hodge structure
$V:=T_{S,\QQ}$ we get a decomposition
$$
T_{S,\QQ}\,\cong\,V_3\otimes W_3'\;\oplus\;V_7\otimes W_7',\qquad
\dim W_3',\;\dim W_7'\in\{2,3\},
$$
where $W_3'$, $W_7'$ are Hodge substructures of  $W_3$, $W_7$ (both of which are $3$ dimensional) respectively with
$\dim (W_i')^{2,0}=1$.

To find the dimension of $W_7'$, we consider the action of $G$ on the subspace $L_S$ of $NS(S)_\QQ$ spanned by classes of lines.
Elements of $G$ permute the lines and using the intersection form one can find an explicit matrix which gives the action of any $g\in G$ on $L_S$.
Thus the representation of $G$ on $L_S$ can be
decomposed into irreducibles.
We found, of course,
that all representations in $H^2(S)_{pr}$ which do not appear in $H^{2,0}$
do appear in $L_{S}$.
The interesting thing is that $V_3$ appears also in $L_{S}$,
hence we can conclude that $\dim W_3'=2$,
but $V_7$ does not appear in $L_{S}$.
Thus we get the decomposition
$$
H^2(S,\QQ)\,=\,V_3\otimes W_3'\,\oplus\,V_7\otimes W_7\,\oplus\,L_S.
$$
Now it only remains to decide whether $\dim W_7'$ is $2$ or $3$ (equivalently, whether $\dim NS(S)_\QQ=202+30$, in this case $L_S\neq NS(S)_\QQ$, or $202$ respectively).
We will see in section \ref{decS} that
$$
T_{S,\QQ}\,\cong\,V_3\otimes W_3'\;\oplus\;V_7\otimes W_7,\qquad
\dim W_3'\,=\,2,\quad \dim W_7\,=\,3,
$$
hence $L_S=NS(S)_\QQ$ and $W_7$ is a simple 3-dimensional Hodge structure.

\section{Quotients of Maschke's surface $S$}\label{quotS}

\subsection{}
In section \ref{dimdh2} we showed that  the transcendental Hodges structure
$T_{S,\QQ}\subset H^2(S,\QQ)$ of $S$ has a summand $V_3\otimes W_3'$ which is isomorphic to $(W_3')^{\oplus 5}$ as a rational Hodge structure.
We will now find desingularizations $\widetilde{U}$, $\widetilde{W}$ of
quotients of $S$ by subgroups of $G$ with transcendental Hodge structure
$T_{\widetilde{U},\QQ}\cong (W_3')^{\oplus 5}$ and $T_{\widetilde{W},\QQ}\cong W_3'$.
We did not succeed in finding a model of a surface with transcendental Hodge structure isomorphic to $W_7'$, the other Hodge substructure in $T_{S,\QQ}$.

\subsection{K3 surfaces and $W_3'$}\label{w3p}
The Heisenberg group $H$ acts through
$\overline{H}:=H/\mu_4\cong (\ZZ/2\ZZ)^4$, an abelian group, on $\PP^3$
and thus on $S$.
This gives a natural splitting
$$
T_{S,\QQ}\,=\,T_{S,\QQ,1}\,\oplus\,
(\oplus_{\chi\neq 1}\,T_{S,\QQ,\chi}),
$$
where the first summand is the subspace of $H$-invariants and
the sum is over the $15$ non-trivial characters of $H$, which factor over
$\overline{H}$.
The character table of $G$ shows that $H$ acts trivially in the $G$-representation $V_3$, so we conclude that $V_3\otimes W_3'\cong
T_{S,\QQ,1}$ and this also identifies $V_7\otimes W_7'\subset T_{S,\QQ}$:
$$
T_{S,\QQ}^H\,=\,T_{S,\QQ,1}\,\cong\,V_3\otimes W_3',\qquad
\oplus_{\chi\neq 1}\,T_{S,\QQ,\chi}\,\cong\,V_7\otimes W_7'.
$$

The desingularization $\widetilde{U}$ of the quotient surface $U:=S/H$ will have
transcendental Hodge substructure $T_{S,\QQ}^H$, hence
$$
T_{\widetilde{U},\QQ}\,\cong\,V_3\otimes W_3'\,\cong\,(W_3')^{\oplus 5}.
$$
A further quotient of $U$ by an involution $\iota_U$ will give us a surface $W$ whose desingularization $\widetilde{W}$ is a K3 surface with $T_{\widetilde{W}}\cong W_3'$.
$$
S\,\longrightarrow\,U\,:=\,S/H\;(\subset\PP^4)\;
{\longrightarrow}\,W\,:=U/\iota_U\;(\subset \PP^3).
$$

\subsection{The surface $U=S/H$}\label{U}
To find a projective model of the quotient surface $U=S/H$ we consider
the subring of $H$-invariant polynomials in $\CC[x_0,x_1,x_2,x_3]$.
It is generated by the following $5$ polynomials of degree $4$:
$$
p_0:\,=\,x_0^4+x_1^4+x_2^4+x_3^4,\quad
p_1\,:=\,2(x_0^2x_1^2+x_2^2x_3^2),\quad
p_2\,:=\,2(x_0^2x_2^2+x_1^2x_3^2),
$$
$$
p_3\,:=\,2(x_0^2x_3^2+x_1^2x_2^2),\quad p_4\,:=\,4x_0x_1x_2x_3,
$$
and we have
$$
\CC[x_0,x_1,x_2,x_3]^H\,=\,\CC[p_0,\ldots,p_4]\,\cong\,\CC[y_0,\ldots,y_4]/(G_I)
$$
where the isomorphism is given by $y_i\mapsto p_i$, with kernel the quartic polynomial
$$
G_I\,:=\, y_4^4\,+\,(y_0^2-y_1^2-y_2^2-y_3^2)y_4^2+y_1^2y_2^2+y_1^2y_3^2+y_2^2y_3^2-2y_0y_1y_2y_3.
$$
Thus the image of the map $\PP^3\rightarrow\PP^4$, $x\mapsto (\ldots:p_i(x):\ldots)$
is the (singular) quotient variety $\PP^3/H$, which is known as the Igusa quartic,
and it is defined by $G_I=0$.

The image $U$ of Maschke's surface $S$ is the intersection of the Igusa quartic with the quadric with defining equation
$$
G_M\,:=\,y_0^2 + 3(y_1^2 + y_2^2 + y_3^2) + 6y_4^2,
$$
so $G_M(p_0,\ldots,p_4)=F$, the defining equation of $S$.
In fact, the group $G/H=S_6$ acts on
$\PP^4=\PP(\CC^5)$ and the representation on $\CC^5$ is irreducible
with unique quadratic invariant $G_M$.
The singular locus of $U:=S/H$ consists of $30$ nodes
which are the images of the fixed points of $H$ in $S$.
By adjunction, we get an isomorphism $H_{\widetilde{U}}^{2,0} \cong H^0(\PP^4,\cO(1))$, the vector space of linear forms on $\PP^4$.

\subsection{The surface $W=U/\iota_U$}\label{W}
We define an involution $\iota_U$ on $U$ and we
give an explicit projective model of the surface $W=U/\iota_U$ as a
12 nodal quartic surface in $\PP^3$.

The matrix $\mbox{diag}(-1,1,1,1)$ lies in $G$ and it induces the map
$\mbox{diag}(1,1,1,1,-1)$ on the basis $p_0,\ldots,p_4$ of the space
of $H$-invariants of  degree $4$, and hence on $\PP^4$.
In any case, it fixes the two defining equations  $G_I$ and $G_M$ of $U$ and hence it induces an involution $\iota_U$ on $U$.

The quotient $W$ of $U$ by $\iota_U$ is the surface in $\PP^3$,
with coordinates $y_0,\ldots,y_3$,
whose defining equation is obtained by substituting
$y_4^2=-(y_0^2 + 3(y_1^2 + y_2^2 + y_3^2)/6$ in $G_I$.
Thus
$$
U\,\longrightarrow\, W\,:=\,U/\langle \mbox{diag}(1,1,1,1,-1)\rangle,\qquad
W\,=\,(H=0)\quad(\subset\PP^3),
$$
where the quartic polynomial $H$ is given by:
$$
H\,:=\,
5y_0^4 + 6y_0^2(y_1^2 + y_2^2 + y_3^2)
-27(y_1^4 +y_2^4+y_3^4)
-90(y_1^2y_2^2 +y_1^2y_3^2 + y_2^2y_3^2)
+ 72y_0y_1y_2y_3.
$$
The automorphism group of $W$ contains a group $G_{24}$ of 24 elements
generated by a subgroup isomorphic to $S_3$, given by permutations of the variables $y_1,y_2,y_3$, and a subgroup isomorphic to $(\ZZ/2\ZZ)^2$,
given by changing an even number of signs of these three variables.
The point $p:=(3:3:\sqrt{-3}:\sqrt{-3})$ is a singular point, a node, of $W$. Its orbit under $G_{24}$ consists of 12 singular points and this set is the singular locus of $W$.

The surface $W$ contains the line $m$ spanned by
$$
m\,:\quad\langle\,(3:0:\sqrt{-3}:0),\;(0:3:0:\sqrt{-3})\,\rangle.
$$
The orbit of this line under $G_{24}$ consists of 12 lines. Another line $m'$
on $W$ is given by
$$
m'\,:\quad\langle\,(\omega-1:-\omega:1:0),\;(-\omega-2:\omega+1:0:1)\,\rangle,
\qquad \omega:=(-1+\sqrt{-3})/2,
$$
so $\omega$ is a primitive cube root of unity.

The desingularization $\widetilde{W}$ of $W$ is a K3 surface.
The subgroup of its N\'eron-Severi group $NS(\widetilde{W})$ generated by the classes of the 12 rational curves over the nodes, the strict transforms of the 12 lines in the $G_{24}$-orbit of $m$ and $m'$ has rank $20$, as an explicit computation of the intersection matrix of these curves shows.
As this is also the maximal value that the rank can have, we conclude that rank$(NS(\widetilde{W}))=20$. This agrees with the fact that
$\dim H^2(\widetilde{W},\QQ)=22$,
$T_{\widetilde{W}}\cong W_3'$ and $\dim W_3'=2$ (cf.\ section \ref{dimdh2}).

In section \ref{LseriesW} we will determine the Galois representation
on $H^2_\et(\widetilde{W},\QQ_\ell)$.

\section{Galois representations and Maschke's surface}\label{H2etS}

\subsection{Decomposing $H_\et^2(S,\QQ_\ell)$}\label{decH2etS}
Maschke's surface $S$ is defined by a polynomial with integer coefficients
and thus there is a representation of the Galois group $G_\QQ:=Gal(\overline{\QQ}/\QQ)$ on the \'etale cohomology groups
$H_\et^2(S,\QQ_\ell)$ (where now $S=S_{\overline{\QQ}}$ is considered as a surface over $\overline{\QQ}$).
The ideal in $\ZZ[x,y,z]$ generated by the partial derivatives of $F(1,x,y,z)$ contains the integer $2^{15}3^25$ (according to Magma).
Hence, using the invariance of $F$ under permutations of the four variables,
the surface $S$ has good reduction at all primes $p$ with $p>5$.

There is a decomposition of $G_\QQ$-representations:
$$
H_\et^2(S,\QQ_\ell)\,=\,T_{S,\ell}\,\oplus\,NS(S)_\ell,
$$
where $NS(S)_\ell$ is the subspace spanned by divisor classes in $S_{\overline{\QQ}}$, so $\dim NS(S)_\ell=\dim NS(S)_\QQ$,
and $T_{S,\ell}=NS(S)_\ell^\perp$.

The action of $G\subset Aut(S)$ is defined by matrices with
entries in $\QQ(i)$. Thus the actions of $G$ and $G_{\QQ(i)}:=Gal(\overline{\QQ}/\QQ(i))\subset G_\QQ$
on $T_{S,\ell}$ commute. Using the results from
section \ref{dimdh2} and the comparison theorem, we find
the splitting, as $G\times G_{\QQ(i)}$-representations:
$$
T_{S,\QQ_\ell}\,\cong\,V_{3,\ell}\otimes W_{3,\ell}'\;\oplus\;
V_{7,\ell}\otimes W_{7,\ell}',\qquad
\dim W_{3,\ell}'\,=\,2,\quad \dim W_7'\in\{2,3\},
$$
we will see that $\dim W_7'=3$ in section \ref{decS}.

As the representations $V_3,V_7$ of $G$ are irreducible and have distinct dimensions,
the two summands are preserved by the $G_\QQ$-action on $T_{S,\QQ_\ell}$.
The $G$-representation on $V_3\otimes W_3'$ factors over $S_6$ and is then
isomorphic to the $S_6$-representation
on $T_{\widetilde{U},\QQ}$.  As $S_6$ acts by matrices with coefficients in $\QQ$ on
$U\subset \PP^4$, we conclude that $W_{3,\ell}'$ is in fact a $G_\QQ$-representation.

Since $W_3'\cong T_{\widetilde{W}}$, as rational Hodge structures, we get by a similar argument that $W_{3,\ell}'$ is the Galois representation on the orthogonal complement of $NS(\widetilde{W})_\ell$ in $H^2_{\et}(\widetilde{W},\QQ_\ell)$. We will describe it
explicitly in section \ref{LseriesW}.

\subsection{The Galois representation on $L_S$}
The $202$ dimensional subspace $L_S$ of $NS(S)_\QQ$ is spanned by the
classes of the lines in $S$. These lines are defined over the field
$\QQ(\zeta_{12},\sqrt{5})$, a subfield of $\QQ(\zeta_{60})$,
where $\zeta_{k}$ is a primitive $k$-th root of unity.
The Galois representation on $L_S$, given by permutations of the classes of the lines, thus factors over $Gal(\QQ(\zeta_{12},\sqrt{5})/\QQ)\cong (\ZZ/2\ZZ)^3$.
This group has $8$ one-dimensional representations, which extend to $G_\QQ$-representations:
$$
\sigma_{a,b,c}\,:=\,\sigma_{1,0,0}^a\sigma_{0,1,0}^b\sigma_{0,0,1}^c\,:\;
G_\QQ\,\longrightarrow\,\{\pm 1\}
$$
for $a,b,c\in\{0,1\}$ and we define the following non-trivial Galois
representations by their kernel, which has index two in $G_\QQ$:
$$
\ker(\sigma_{1,0,0})\,=\,G_{\QQ(\sqrt{-1})},\qquad
\ker(\sigma_{0,1,0})\,=\,G_{\QQ(\sqrt{-3})},\qquad
\ker(\sigma_{0,0,1})\,=\,G_{\QQ(\sqrt{5})}.
$$
Thus for a Frobenius element
$F_p\in G_\QQ$ at the prime $p$ one has:
$\sigma_{1,0,0}(F_p)=-1$ iff $p \equiv 3$ mod 4,
$\sigma_{0,1,0}(F_p)=-1$ iff $p \equiv 2$ mod 3 and
$\sigma_{0,0,1}(F_p)=-1$ iff $p \equiv 2,3$ mod 5.

With a computer we found that
$$
L_S\,\cong\,
\sigma_{0,0,0}^{44}\oplus\sigma_{0,0,1}^{28}\oplus
\sigma_{0,1,0}^{28}\oplus\sigma_{1,0,0}^{42}\oplus
\sigma_{1,0,1}^{33}\oplus \sigma_{1,1,0}^{27},
$$
note that $\dim L_S=202=44+28+28+42+33+27$.
The $G_\QQ$-representation on $L_{S,\ell}\subset H^2_\et(S,\QQ_\ell)$ is
then the tensor product of the representation on $L_S$ above
with  $\QQ_\ell(-1)$. In particular, for $p>5$ a Frobenius element
$F_p\in G_\QQ$ at the prime $p$ will act as $\sigma_{a,b,c}(F_p)p$ on the summand $\sigma_{a,b,c}\otimes \QQ_\ell(-1)$.
This allows one to compute the trace of $F_p$ on $L_{S,\ell}$ explicitly.

\subsection{The Galois representation on $H^2_{\et}(\widetilde{W},\QQ_\ell)$.} \label{LseriesW}
The surface $\widetilde{W}$, defined in section \ref{W}, has good reduction at primes $p>5$.
The Galois representation on $H^2_{\et}(\widetilde{W},\QQ_\ell)$
splits into the 20-dimensional representation on
$NS(\widetilde{W})_\ell$
and a two dimensional representation on the orthogonal complement
$T_{\widetilde{W},\ell}\cong W_{3,\ell}'$.
Using the explicit description of the generators of $NS(\widetilde{W})_\QQ$
one finds that all generators are rational over $\QQ(\sqrt{-3})$ and that
the non-trivial element in the Galois group $Gal(\QQ(\sqrt{-3})/\QQ)$
has $10$ eigenvalues $+1$ and $10$ eigenvalues $-1$ on this sublattice.
Hence
$$
NS(\widetilde{W})_\ell\,\cong\,
\QQ(-1)^{10}\,\oplus\;\;\sigma_{0,1,0}\otimes\QQ_\ell(-1)^{10}.
$$
Moreover, we found that
the determinant of the intersection matrix on a basis of sublattice is,
up to a square, equal to $-15$.
The theory of the arithmetic of singular K3 surfaces
then asserts that the representation on $T_{\widetilde{W}_\ell}$
is determined by a
Hecke character of the imaginary quadratic field $L:=\QQ(\sqrt{-15})$.

To determine the Hecke character,
a computer counted the number of points $\sharp W(\FF_p)$
on $W$ in $\PP^3(\FF_q)$ for many small prime powers $q$.
If $q\equiv 1$ mod 3, the 12 nodes are rational over $\FF_q$ and thus $\sharp \widetilde{W}(\FF_q)=\sharp W(\FF_q)+12p$, else
$\sharp \widetilde{W}(\FF_q)=\sharp W(\FF_q)$.
Then Lefschetz's fixed point formula implies:
$$
\sharp \widetilde{W}(\FF_q)\,=\,1+10(1+\sigma_{0,1,0}(F_q))q+a_q+q^2,
\qquad
a_q:=Trace(F_p^k|T_{\widetilde{W},\ell}),\qquad \mbox{where}\quad q=p^k.
$$
The values of various $a_p$ are listed below.

Let $\cO_L:=\ZZ[\alpha]$, $\alpha:=(1+\sqrt{-15})/2$
be the ring of integers of  $L$.
Its class number is $2$, hence the square of any of its ideals is principal.
As the units in $\cO_L$ are $\pm 1$,
a generator of a principal ideal is unique up to sign.
As the minimal polynomial of $\alpha$ is $f:=x^2-x+4$, and $f\equiv x^2+2x+1\equiv (x+1)^2$ mod 3, there is a surjective homomorphism
$$
\phi_3:\,\cO_L:=\ZZ[\alpha]\,\longrightarrow\,\ZZ/3\ZZ,\qquad a+b\alpha\,\longmapsto\,
a-b\;\mbox{mod}\;3
$$
with $\phi_3(-1)=-1$. Thus any principal ideal $I=\beta\cO_L$, prime with $3\cO_L$, has a unique generator $\beta$ such that $\phi_3(\beta)=+1$. We define a Hecke character $\chi$ on $\cO_L$ by defining
$$
\chi(I)\,=\,\beta,\qquad\mbox{where}\quad I^2=(\beta),\quad \phi_3(\beta)=1.
$$

To be explicit, if $p=\gop\overline{\gop}$ is a prime which splits in $\cO_L$,
then $\gop^2=(a+b\alpha)$ for some integers $a,b$, with $a,b\neq 0$,
which can be found by considering the norm:
$$
p^2\,=\,N(a+b\alpha)\,:=\,(a+b\alpha)\overline{(a+b\alpha)}\,=\,a^2-ab+4b^2.
$$
Changing signs, if necessary, to obtain $a-b\equiv 1$ mod 3, one thus has
$\chi(\gop)=a+b\alpha$ and $\chi(\gop)+\overline{\chi(\gop)}=2a+b$.
For example, if $p=17$ and $\wp_{17}$ is one of the primes over $17$, then
$$
11^2-11\cdot 8+8^2=17^2,\quad
\phi_3(11-8\alpha)\,=\,19\,\equiv\, 1\, \mbox{mod}\, 3,
\quad\mbox{so}\quad
\chi(\gop_{17})\,=\,11-8\alpha.
$$

This Hecke character determines a one dimensional
$\ell$-adic Galois representation of $G_L:=Gal(\overline{\QQ}/L)$
which induces a two dimensional
Galois $\ell$-adic representation of $G_\QQ$.
This is the representation of $G_\QQ$ on
$T_{\widetilde{W}_\ell}=W_{3,\ell}'$ which we found.
In particular, the trace $a_p$  of $F_p$ in this representation is:
$$
a_p\,=\,\left\{\begin{array}{rl} \chi(\gop)+\overline{\chi(\gop)}\qquad &\mbox{if $p$ splits}\\
0\qquad&\mbox{if $p\cO_L$ is prime}.\end{array}\right.
$$
For example, $\chi(\gop_{17})\,=\,11-8\alpha$ so
we get $a_{17}\,=\,2\cdot11+(-8)=14$, below are some other values of $a_p$.

\

{\renewcommand{\arraystretch}{1.2}
$$
\begin{array}{|r|r|r|r|r|r|r|r|r|r|r|r|r|r|}\hline
p&7&11&13&17&19&23&29&31&\ldots&79&83&89&97 \\
\hline
a_p&0&0&0&14&-22&-34&0&2&\ldots&98&-154&0&0\\
\hline
\end{array}
$$
}

\

The $a_p$ are also the Fourier coefficients of
an elliptic modular form $f$ of weight three (cf.\ \cite{S}).
In this case $f$ is one of the two newforms
with Dirichlet character $p\mapsto \left(\frac{-15}{p}\right)$, it has level $N=15$.

Before continuing with the study of $H^2_\et(S,\QQ_\ell)$ in section \ref{decS},
we make a digression to discuss various surfaces related to $S$.

\subsection{Remark} \label{bS}
The polynomial $F$ defining $S$ is a quartic in the $x_i^2$,
hence $S$ is a ramified cover of the surface $\overline{S}$
defined by
$$
\overline{S}:\qquad
\sum_{i=0}^3 x_i^4\,+\,14\sum_{i<j}x_i^2x_j^2\,+\,168x_0x_1x_2x_3.
$$
The surface $\overline{S}$ is smooth and hence it is a K3 surface.
It is the quotient of $S$ by the subgroup $K\cong (\ZZ/2\ZZ)^4$ of diagonal matrices with entries $\pm 1$ of $G$ and thus
$H^2(\overline{S},\QQ)\cong H^2(S,\QQ)^H$.

The Galois representation
on the transcendental lattice $T_{\overline{S},\ell}$ is the same as
the one on $T_{\widetilde{W},\ell}$:
$$
\sharp \overline{S}(\FF_p)\,=\,1+n_pp+a_p+p^2,
$$
where $n_pp$ is the trace of the Frobenius $F_p$ on the N\'eron-Severi group
of $\overline{S}$ , with
$$
n_p\,:=\,(4\sigma_{0,0,0}\oplus4\sigma_{0,1,0}\oplus
6\sigma_{1,0,0}\oplus 3\sigma_{1,0,1}\oplus 3\sigma_{1,1,0})(F_p).
$$

\subsection{The Galois representation on $H^2_{\et}(\widetilde{U},\QQ_\ell)$}\label{LU}
Using the action of $S_6$ on $H^2(\widetilde{U},\QQ)$, we already concluded that $T_{\widetilde{U}}\cong W_{3}'^5$.
As the automorphisms of $\widetilde{U}$ corresponding to the elements of $S_6$ are defined over $\QQ$ this implies that the $G_\QQ$-representation on
$H^2_{\et}(\widetilde{U},\QQ_\ell)$ splits into $T_{\widetilde{U},\ell}\cong T_{\widetilde{W},\ell}^5$ and
a direct sum of $\QQ_\ell(-1)$'s twisted by Dirichlet characters, unramified for primes $p>5$.

As $U$ is a complete intersection of type $(2,4)$ in $\PP^4$
with only nodes as singularities,
the Euler characteristic $\chi(\widetilde{U})$ of
$\widetilde{U}$ is the same as the one of a smooth
complete intersection of type $(2,4)$, so $\chi(\widetilde{U})=64$.
As the cohomology in odd degree is zero, we get $\dim H^2(\widetilde{U},\ZZ)=62$ and thus $\dim NS(\widetilde{U})_\QQ=52$.

Note that $U$ has 30 nodes, defined over $\QQ(\sqrt{-3})$, hence
$\sharp \widetilde{U}(\FF_p)
=\sharp U(\FF_p)+30p$ if $p\equiv 1$ mod 3 and else $\sharp \widetilde{U}(\FF_p)
=\sharp U(\FF_p)$, for $p>5$.
Using a computer and the Lefschetz
fixed point formula one finds, for many small primes $p$ with $p>5$,
$$
\sharp \widetilde{U}(\FF_p)\,=\,
1+(26+25\sigma_{0,1,0}(F_p)+\sigma_{0,0,1}(F_p))p+5a_p+p^2,
$$
hence this determines the $G_\QQ$-representation on $H^2_{\et}(\widetilde{U},\QQ_\ell)$.

\

\subsection{The Galois representation $W_{7,\ell}$}\label{decS}
From sections \ref{dimdh2}, \ref{decH2etS} we have the following decomposition of $H^2_{\et}(S,\QQ_\ell)$ as
$G\times G_{\QQ(i)}$-representations:
$$
H^2_{\et}(S,\QQ_\ell)\,=\,V_{3,\ell}\otimes W_{3,\ell}'\,\oplus \,
V_{7,\ell}\otimes W_{7,\ell}\,\oplus\,L_{S,\ell}
$$
where $W_{7,\ell}$ is three dimensional. As $\dim V_{3,\ell}=5$, we have
 $Tr(F_q|V_{3,\ell}\otimes W_{3,\ell}')=5a_q$ with
$a_q:=Tr(F_q|T_{\widetilde{W},\ell})$ as in section \ref{LseriesW}.
Using that $\dim V_{7,\ell}=30$ and taking prime numbers
 $q\equiv 1$ mod 4,
the Lefschetz fixed point formula gives:
$$
\sharp S(\FF_q)\,=\,
1+5a_q+30b_q+tr(F_q|L_{S,\ell})+q^2,
\qquad b_q:=Tr(F_q|W_{7,\ell}).
$$
Some $b_q$, for primes $q$, are listed in section \ref{LsS} below.

The $G_{\QQ(i)}$-representation on $W_{7,\ell}$ is compatible with the non-degenerate quadratic form induced by the intersection form on
$H^2_\et(S,\QQ_\ell)$.
If $W_{7,\ell}$ is reducible, then it must be the direct sum of
a one and a two dimensional representation,
which is compatible with a non-degenerate quadratic form.
The orthogonal group $O(2)$ is essentially abelian, and
therefore the two dimensional representation is of CM type, in particular the eigenvalues of the Frobenius elements must be in a fixed CM field
of degree at most two over $\QQ(i)$.
However, we checked that this is not the case,
using that the characteristic polynomial
for $F_p$ with $p$ a prime which is $1$ mod 4,
is given by $x^3-b_px^2+\,\epsilon_pb_pp  x-\epsilon_pp^3$ with
$\epsilon_p\in\{\pm 1\}$ (cf.\ section \ref{LsS}).

Therefore the $G_{\QQ(i)}$-representation $W_{7,\ell}$ is irreducible
and hence $L_S=NS(S)_\QQ$ as $\QQ$-vector spaces.
This again implies that $W_7$ is a simple Hodge structure of dimension three.

In the recent paper \cite{RS} another approach to solve a similar problem is given.

\subsection{The Galois representation on $H^2_{\et}(S,\QQ_\ell)$.}\label{LsS}
We conclude our study of the arithmetic of $S$ by determining the
$G_\QQ$-representation on $H^2_{\et}(S,\QQ_\ell)$ in terms of known
representations and a three dimensional $G_\QQ$-representation which
restricts to the $G_{\QQ(i)}$-representation $W_{7,\ell}$.

Upon restriction to the Heisenberg group $H$,
the irreducible 30 dimensional $G$-representation $V_7$
splits into a direct sum
$$
V_7\,\cong\,\oplus_{\chi\neq 1} V_{7,\chi},\qquad \dim V_{7,\chi}\,=\,2,
$$
where the sum is over the 15 non-trivial one dimensional representations of $H$.
As $H$ is defined over $\QQ$, we get a $H\times G_\QQ$-subrepresentation
$$
\oplus_{\chi\neq 1}\;
V_{7,\chi,\ell}\otimes W_{7,\chi,\ell}\;\subset\;H^2_{\et}(S,\QQ_\ell),
$$
where the $W_{7,\chi,\ell}$ are three dimensional $G_\QQ$-representations.

The group $G$ acts on $H$ and it permutes these 15 one dimensional representations transitively. Thus the $G_\QQ$-representations $W_{7,\chi,\ell}$ are all isomorphic when restricted to the subgroup $G_{\QQ(i)}$,
and in fact they are all isomorphic to $W_{7,\ell}$.
Fix one non-trivial $\chi_0$, then given a non-trivial $\chi$ we have that
either
$$
W_{7,\chi,\ell}\,\cong\,W_{7,\chi_0,\ell}\quad\mbox{or}\quad
W_{7,\chi,\ell}\,\cong\,W_{7,\chi_0,\ell}\otimes\sigma_{1,0,0}.
$$
As $\dim V_{7,\chi,\ell}=2$, we get a decomposition, as $G_\QQ$-representations:
$$
H^2_{\et}(S,\QQ_\ell)\,=\, (W_{3,\ell}')^{\oplus 5}\,\oplus \,
(W_{7,\chi_0,\ell})^{\oplus 2a}\,\oplus \,
(W_{7,\chi_0,\ell}\otimes\sigma_{1,0,0})^{\oplus 2b} \,\oplus\, L_{S,\ell},\qquad a+b=15,
$$
for some integers $a,b$ and a certain $3$-dimensional $G_\QQ$-representation
$W_{7,\chi_0,\ell}$. Replacing $W_{7,\chi_0,\ell}$ by $W_{7,\chi_0,\ell}\otimes\sigma_{1,0,0}$ if necessary, we may assume that $a>b$.
This determines the $G_\QQ$-representation $W_{7,\chi_0,\ell}$ uniquely.

To determine $a,b$, we use the Lefschetz fixed point formula:
$$
\sharp S(\FF_q)\,=\,
1+5a_q+(2a+2b\sigma_{1,0,0}(F_q))b_p+tr(F_q,L_{S,\ell})+q^2,
\qquad b_q:=Tr(F_q|W_{7,\chi_0,\ell}).
$$
In case $q\equiv 1$ mod 4 one has $\sigma_{1,0,0}(F_q)=1$,
so we recover the formula from section \ref{decS}.
In case $q\equiv 3$ mod 4, $\sigma_{1,0,0}(F_q)=-1$
and the computer finds $2(a-b)b_q$.
For small primes $q$, we found that the g.c.d.\ of these integers is $6$.
Hence
$$
a-b=1
\qquad\mbox{or}\qquad
a-b=3.
$$

We now exclude the case $a-b=1$.
For $q=11$ we found $2(a-b)b_{11}=-78$, hence $b_{11} =-39$ if $a-b=1$.
But $|b_{11}|\leq 3q=33$ since $b_{11}$ is the sum of the three eigenvalues of $F_q$ on $W_{7,\chi_0,\ell}\subset H^2_{\et}(S,\QQ_\ell)$.
Therefore we must have $a-b=3$ and thus
$$
a=9,\qquad b=6.
$$
The following is a table of some of the $b_p$ for primes $p$:

{\renewcommand{\arraystretch}{1.2}
$$
\begin{array}{|r|r|r|r|r|r|r|r|r|r|r|r|r|r|}\hline
p&7&11&13&17&19&23&29&31&37&41&43&47&53 \\
\hline
b_p&-7&-13&-11&5&7&13&-21&19&13&-9&29&-11&-55\\
\hline
\end{array}
$$
}

\

To determine the characteristic polynomial of $F_p$ on $W_{7,\chi_0,\ell}$,
knowing only $b_p$, we proceed as follows.
First of all, if $\alpha$ is an eigenvalue of $F_q$ on $W_{7,\chi_0,\ell}$
then so is $q^2/\alpha=\bar{\alpha}$. As $\dim W'_{7,\chi_0,\ell}=3$ this implies that one of the eigenvalues is $\pm q$.
Let $\alpha,\bar{\alpha},\epsilon_pp$ be the eigenvalues of $F_p$, with
$\epsilon_p\in\{\pm 1\}$, then $\det(F_p)=\epsilon_p p^3$.
As the determinant is also a $G_\QQ$-representation, unramified for $p>5$,
this allows us to determine it explicitly once we know it for small primes.

To find $\epsilon_p$ for small primes we used the following identity involving
$b_p,b_{p^2}$, which we determine by counting points, and $\epsilon_p$:
$$
b_p\,=\,\alpha+\bar{\alpha}+\epsilon_p p,\quad
b_{p^2}\,=\,\alpha^2+\bar{\alpha}^2+ p^2\,=\,
b_p^2-2p(b_p-\epsilon_p p)-2p^2.
$$
We found that
$$
\epsilon_p\,=\,\sigma_{1,0,1}(F_p),
$$
hence the determinant of the $G_\QQ$-representation on $W_{7,\chi_0,\ell}$
is $\sigma_{1,0,1}$.
The characteristic polynomial $f_p$ of $F_p$ on $W_{7,\chi_0,\ell}$ is thus determined
by $b_p$:
$$
f_p\,=\,x^3-b_px^2+\,\epsilon_pb_pp  x-\epsilon_pp^3\,=\,(x^2-(b_p-\epsilon_pp)x  +p^2)(x-\epsilon_p p).
$$

\subsection{Remark}\label{sym2}
M.\ Sch\"utt pointed out that there is a K3 surface ${\mathcal X}_{3,2/3}$,
defined over $\QQ$, whose $G_\QQ$-representation on $H^2_\et$ has a
rank three summand $T_{3,2/3,\ell}$ which has the same characteristic polynomial of $F_p$ for
primes $p$ with $7\leq p\leq 131$ as $W_{7,\chi_0,\ell}$. This surface has an elliptic fibration
given by:
$$
{\mathcal X}_{3,r}:\qquad
y^2 = x^3- t^2 (r^2 t-1-2 r) x^2-2 (t+1) t^3 r (r t-1) x-(t+1)^2 t^5 r^2,
$$
with $r=2/3$.
This family of elliptic surfaces is studied in \cite{GS}, section 6, where it is shown that the ${\mathcal X}_{3,r}$ are birationally isomorphic to quotients of products of two isogeneous elliptic curves.  For $r=2/3$ there is an elliptic curve  over $\QQ(\sqrt{-5}, \sqrt{-15})$, that is 2-isogeneous to its $Gal(\QQ(\sqrt{-5})/\QQ)$
conjugate and 3-isogeneous to its $Gal(\QQ(\sqrt{-15})/\QQ)$ conjugate, which produces the K3 surface. It is the curve $C^{(a)}$ in the family parametrised
by $X^\ast(6)$ with parameter value $a=-16/5$ in
\cite{Q}, p.312.

If there is indeed an isomorphism of $G_\QQ$-representations
$T_{3,2/3,\ell}\cong W_{7,\chi_0,\ell}$, then
the Tate conjecture predicts a correspondence, defined over $\QQ$,
between $S$ and ${\mathcal X}_{3,2/3}$.

Using the representation of $G$ on $H^{2,0}(S)$ it is not hard to find a subgroup $K\subset G$ such that $\dim H^{2,0}(S)^K=1$
and $H^{2,0}(S)^K\subset  (V_7\otimes W_7)\otimes\CC$. Thus any desingularization of the quotient surface $S/K$ has transcendental Hodge structure isomorphic to the rational Hodge structure $W_7$, but we have not been able to find an explicit model of this quotient surface.

\

\section{Maschke's double octic}\label{mdo}

\subsection{The cohomology of $X$}
Maschke's double octic is the Calabi-Yau threefold which is the double cover
$$
\pi:\,X\,\longrightarrow\,\PP^3
$$
branched along Maschke's (smooth) degree $8$ surface $S\subset\PP^3$.
As the surface $S$ has Euler-Poincar\'e characteristic $\chi_{top}=304$,
$X$ has Euler-Poincar\'e characteristic
$$
\chi_{top}(X)\,=\,2\cdot \chi_{top}(\PP^3)\,-\,\chi_{top}(S)=8-304=-296.
$$
Results of Lazarsfeld (\cite{Laz}, Thm.\ 1, Prop.\ 3.1) imply that
$$
H^i(\PP^3,\ZZ)\,\stackrel{\cong}{\longrightarrow}\,H^i(X,\ZZ),\qquad i=0,1,2.
$$
Recall that $H^i(\PP^3,\ZZ)=0$ if $i=1$ and $H^i(\PP^3,\ZZ)=\ZZ$ if $i=0,2$,
and that $X$ is a Calabi-Yau threefold, hence we get:
$$
h^3(X)\,=\,300,\qquad h^{3,0}(X)\,=\,1,\quad h^{2,1}(X)\,=\,149.
$$
As $X$ is a hypersurface in
weighted projective space $WP(1,1,1,1,4)$, one can also use toric methods (cf.\ \cite{CK}, 4.1.3)
to compute the $h^i$.

\subsection{The action of $G$ on $H^3(X,\QQ)$}\label{decXG}
Using our generalization of Ch\^enevert's formula, cf.\ Proposition \ref{propch},
we found that $H^3(X,\QQ)$ is the direct sum of $7$ distinct irreducible $G$-representations, each with multiplicity two and of
$$
\mbox{dimension}\qquad 1,\quad 5,\quad 9,\quad 15,
\quad 30,\quad 45,\quad 45.
$$
The characters of these representations are integer valued,
thus, labeling these irreducible $G$-representations by their dimension,
we have (cf.\ section \ref{cga})
$$
H^3(X,\QQ)\,\cong\,
V_1\otimes W_1\;\oplus\;
V_5\otimes W_5\;\oplus\;
V_9\otimes W_9\;\oplus\;
V_{15}\otimes W_{15}\;\oplus\;
V_{30}\otimes W_{30}\;\oplus\;
V_{45}\otimes W_{45}\;\oplus\;
V_{45}'\otimes W_{45}',
$$
where the $W_i$ are rational Hodge structures of dimension two.
Obviously, $H^{3,0}(X)$ and $H^{0,3}(X)$ are one dimensional subrepresentations in $H^3(X,\CC)$ and the $G$-action on $V_1$ factors over the sign representation of $S_6$, hence $W_1^{3,0}\cong H^{3,0}(X)$.
The other six $W_i$ thus have $W_i^{3,0}=0$, $\dim W_i^{2,1}=1$, so they are Tate twists of Hodge structures of weight one.

The rational Hodge structure $H^3(X,\QQ)$ is thus isomorphic
to $\oplus W_i^{m_i}$, where $m_i=\dim V_i$,
and all seven $W_i$ have dimension two.
We define (the isogeny class of) the elliptic curve  $E_i$
by $E_i=W_i^{p,q}/\Lambda_i$, where $p<q$,
$\Lambda_i\subset W_i$, $\Lambda_i\cong\ZZ^2$ and $\Lambda_i\otimes_\ZZ\QQ=W_i$.
There are natural isomorphisms of $\QQ$-vector spaces
$$
H^1(E_i,\QQ)\,\longrightarrow\, W_i,\qquad (E_i=W_i^{p,q}/\Lambda_i),
$$
which, for $i\neq 1$, are morphisms of rational Hodge structures
with
$H^{p,q}(E_i)\stackrel{\cong}{\rightarrow}W_i^{p+1,q+1}$.
For $i=1$ there cannot exist a morphism of Hodge structures
between $H^1(E_1,\QQ)$ and $W_1$ because $W_1^{2,1}=0$.
In particular,  the Griffiths intermediate Jacobian $J(X)$ of $X$,
a complex torus associated to $H^3(X,\QQ)$ (cf.\ section \ref{ajm}),
is isogeneous to a product $\prod E_i^{m_i}$ of elliptic curves.

\subsection{The action of $H$ on $H^3(X,\QQ)$}\label{decXH}
The Heisenberg group $H\subset G$ acts through its abelian quotient
$\overline{H}\cong (\ZZ/2\ZZ)^4$
on $X$ and thus on the cohomology groups $H^i(X,\QQ)$. We found that
$$
H^3(X,\QQ)\,=\,H^3(X,\QQ)^{H}
\,\oplus\,\left(
\oplus_{\chi\neq 1}
H^3(X,\QQ)_\chi\,\right),
$$
where the sum is over the $15$ non-trivial characters of $H$
and
$$
\dim H^3(X,\QQ)^{H}\,=\,30,\quad
\dim H^3(X,\QQ)_\chi\,=\,18.
$$
Comparing this to the $G$-decomposition above, we find that the $H$-invariants are:
$$
H^3(X,\QQ)^{{H}}\,\cong\,
V_1\otimes W_1\;\,\oplus\;V_5\otimes W_5\;\oplus\;
V_9\otimes W_9.
$$
In the next section we will consider the Heisenberg quotient $Y:=X/H$ of $X$
and its desingularization $\widetilde{Y}$.

\section{The arithmetic of $X$}

\subsection{The Heisenberg quotient $Y$ of $X$}\label{Y}
The CY threefold $X$ is the double cover of $\PP^3$ branched along Maschke's
surface $S$.
The Heisenberg group $H$ acts on $X$, through its action on $\PP^3$.
The quotient of $\PP^3$ by $H$ is the Igusa quartic threefold $Z$ in $\PP^4$, cf.\ section \ref{U}, so there is a double cover
$$
Y\,=\,X/H\,\longrightarrow\,Z\,\cong\,\PP^3/H,
$$
which is branched along the image $U$ of $S$ in $Z$,
which is the intersection of  $Z$ with a quadric in $\PP^4$.

The variety $Z$ is singular along $15$ lines, defined over $\QQ$,
which are the images of the fixed lines in $\PP^3$ of elements in $H$.
These $15$ lines in $\mbox{Sing}(Z)$ intersect, three at the time, in $15$ points.
A desingularization $\hZ$ is obtained by blowing up the singular locus (cf.\ \cite{LW}).
The fiber of $\hZ\rightarrow Z$ over a point $p\in \mbox{Sing}(Z)$
which lies on exactly one line is a $\PP^1$, but the fiber over a point
in three lines consists of three $\PP^1$'s meeting in one point.

The fixed point set in $\PP^3$ of a non-trivial element $h\in H$
is the union of two lines and the stabilizer in $H$ of a general point in such a line is the subgroup generated by $h$. Locally the action of $h$ is  given by the action of $diag(1,-1,-1)$ on $\CC^3$, which has quotient $\CC\times A_1$ where $A_1\cong Spec(\CC[u,v,w]/(uw-v^2))$. This gives the local description of the singular locus of $Z=\PP^3/H$.
It follows, using adjunction, that the strict transform in $\hZ$ of a (linear) hyperplane section of $Z$ is an anti-canonical divisor in $\hZ$.

The variety $\hY$ is the double cover of $\hZ$ branched along the strict transform of the image $U$ of $S$ in $Z$, which is isomorphic to $\widetilde{U}$. As the smooth surface $\widetilde{U}$ is a divisor in $|-2K_{\hZ}|$, the threefold $\hY$ has trivial canonical bundle.

The variety $\hY$ is birationally isomorphic to the crepant resolution $\tilde{Y}'$ of $X/H$, which is a Calabi-Yau threefold.
Thus the Hodge numbers of $\hY$ and $\hY'$ are the same (see \cite{vv})
and we determined them using orbifold cohomology (see \cite{CR}):
$$
h^1(\hY)=0,\quad h^2(\hY)\,=\,h^{1,1}(\hY)\,=16,\qquad
h^3(\hY)\,=\,30,\quad h^{3,0}(\hY)=1,\quad h^{2,1}(\hY)=14.
$$
Comparing this with the $H$-invariants in $H^3(X,\QQ_\ell)$,
we conclude that
$$
H^3_{\et}(\hY,\QQ_\ell)\,\cong\, H^3_{\et}(X,\QQ_\ell)^H\,\cong\,
V_{1,\ell}\otimes W_{1,\ell}\;\,\oplus\;V_{5,\ell}\otimes W_{5,\ell}\;\oplus\;
V_{9,\ell}\otimes W_{9,\ell},
$$
where the $W_{i,\ell}$ are two dimensional $\QQ_\ell$ vector spaces.

As the $S_6$-representation on $H_\et^3(\hY,\QQ_\ell)$
is induced by the action of $S_6$ on $\PP^4$,
where it acts by matrices with rational coefficients, the
$W_{i,\ell}$ are two dimensional $G_\QQ$-representations:
$$
\sigma_i:\,G_\QQ\,:=\,Gal(\overline{\QQ}/\QQ)\,
\longrightarrow\, GL(W_{i,\ell})\,\cong\,GL(2,\QQ_\ell)
$$
for $i=1,5,9$.

Recent modularity results (\cite{GY},\cite{K},\cite{D})
imply that the Galois representation $\sigma_1$
corresponds to a newform of weight $4$ on $\Gamma_0(N)$ for an integer $N$
which is divisible only by primes where $\hY$ has bad reduction.
The representations $\sigma_5,\sigma_9$
are expected to be Tate twists of Galois representations associated
to elliptic curves defined over $\QQ$,  these curves should become isogeneous,
over $\CC$, to the curves $E_5$, $E_9$ from section \ref{decXG}.
Hence, by Wiles' theorem,
these Galois representations  should correspond to newforms of weight two on $\Gamma_0(N_i)$ ($i=5,9$)
for certain integers $N_i$
divisible only by primes where $\hY$ has bad reduction.

\subsection{The Galois representation on $H^3_{\et}(\hY,\QQ_\ell)$}\label{comps}
From  the description of $\hZ$ as the blow up of $Z$,
it follows in particular that $H^2_\et(\hZ,\QQ_\ell)$ is generated by classes of
divisors defined over $\QQ$, hence the Galois representation on $H^2_\et$
(and by duality, also the one on $H^4_\et$) are direct sums of the
Tate representations $\QQ_\ell(-1)$ (and $\QQ_\ell(-2)$ respectively),
where, as usual, $\QQ_\ell(-k)$ is the one dimensional $Gal(\overline{\QQ}/\QQ)$
representation on which the Frobenius $F_q$ acts as multiplication by $q^k$.
In particular:
$$
\sharp \hZ(\FF_q)\,=\,\sharp Z(\FF_q)\,+\,15(q^2+q).
$$

Using the 2:1 map $\hY\rightarrow \hZ$ it then follows that also
$$
\sharp \hY(\FF_q)\,=\,\sharp Y(\FF_q)\,+\,15(q^2+q).
$$
(Use that over each singular line ($\cong\PP^1$ in $Z$) there is a $\PP^1$-bundle over $\PP^1$ in $\hat{Z}$ and after taking the double cover, its preimage is a
$\PP^1$-bundle over a $\PP^1$ in $\hY$.)

From the Lefschetz fixed point formula one has:
$$
\sharp \hY(\FF_q)\,=\,\sum_{i=0}^6\,(-1)^i tr(F_q|H^i_\et(\hY,\QQ_\ell)).
$$
The double cover $\hY\rightarrow \hZ$ induces an isomorphism on $H^2$ and $H^4$
and this map is defined over $\QQ$, so we get:
$$
\sum_{i=0}^6\,(-1)^i tr(F_q|H^i_\et(\hY,\QQ_\ell))\,=\,1+16(q+q^2)+q^3-
tr(F_q|H^3_\et(\hY,\QQ_\ell)).
$$
The trace of the Frobenius $F_q$ on $H^3_\et(\hY)$ can thus be determined by counting points on the (singular) variety $Y\subset\PP^5$:
$$
tr(F_q|H^3_\et(\hY,\QQ_\ell))\,=\,-\sharp Y(\FF_q)\,+\,1+q+q^2+q^3.
$$
Computer computations lead to the following table
{\renewcommand{\arraystretch}{1.2}
$$
\begin{array}{|r|r|r|r|r|r|r|r||r|r|r|r|}\hline
q&7&11&13&17&19&23&29&7^2&\ldots&19^2 \\
\hline
tr(F_q|H^3)
&0&180&210&-90&-1020&-1560&1410&-10290&\ldots&-122970\\
\hline
\end{array}
$$
}

The table below lists the Fourier coefficients of some newforms on $\Gamma_0(N)$
of weight $k$, the names of the weight two forms are those from Table 3 in \cite{modfun4}.

{\renewcommand{\arraystretch}{1.2}
$$
\begin{array}{|r|r|r|r|r|r|r|r|r|r|r|r|r|r|r|}\hline
&p&7&11&13&17&19&23&29&31&\ldots&79&83&89&97 \\
\hline
f120, k=4&a_p&0&4&54&114&44&96&134&-272&\ldots
&688&1188&-694&-1726 \\
\hline
f24B,k=2 &b_p&0 &4&-2&2&-4&-8&6&8&\ldots
&-8&-4&-6&2\\
\hline
f120E,k=2&c_p&0&-4&6&-6&-4&0&-2&-8&\ldots
&-8&-12&10&2\\
\hline
f15C,k=2&d_p&0&-4&-2&2&4&0&-2&0&\ldots
&0&12&-6&2\\
\hline
\end{array}
$$
}

\

Then one can verify that for $7\leq p\leq 97$ one also has:
$$
tr(F_p|H^3_\et(\hY,\QQ_\ell)) \,=\, a_p\,+\,p(\,9b_p\,+\,5c_p).
$$
This leads us to conjecture that
$$
W_{1,\ell}\,\stackrel{?}{\cong}\, V_{f{120},\ell},\qquad
W_{5,\ell}\,\stackrel{?}{\cong}\, V_{f{120E,\ell}}(-1),\qquad
W_{9,\ell}\,\stackrel{?}{\cong}\, V_{f{24B,\ell}}(-1),
$$
where $V_{g,\ell}$ denotes the $\ell$-adic Galois representation associated to the newform $g$.

To find the $a_p$, we assumed that $\sigma_5$ and $\sigma_9$
are Tate twists of Galois representations,
in particular that $tr(F_p|W_{i,\ell})$ is a multiple of $p$ for $i=5,9$. Thus, by counting points, one can find $a_p$ mod $p$.
Comparing with the Fourier coefficients of weight four newforms of level $N=2^a3^b5^c$ for small $a,b,c$ (using Magma) we found that those of $f120$ match perfectly.

Assuming that the $a_p$ are determined correctly, one can determine
$b_p,c_p$ by counting points over $\FF_p$ and $\FF_{p^2}$.
In fact, let $\beta,\bar{\beta},\gamma,\bar{\gamma}\in\CC$ be the eigenvalues of the Frobenius $F_p$ in the two-dimensional weight two representations, so
$$
T^2-b_pT+p\,=\,(T-\beta)(T-\bar{\beta}),\qquad
T^2-c_pT+p\,=\,(T-\gamma)(T-\bar{\gamma}).
$$
Then the traces of $F_p^2$ are given by:
$$
b_{p^2}\,=\,\beta^2+\bar{\beta}^2\,=\,b_p^2-2p,\qquad
c_{p^2}\,=\,\gamma^2+\bar{\gamma}^2\,=\,c_p^2-2p
$$
which allows one to find a degree two polynomial whose zeroes are
$b_p$ and $c_p$ respectively.
The polynomials turn out to have a double zero or a unique zero which is an integer for the case $p=7,11,13,17,19$ thus allowing one to determine $b_p,c_p$ for these primes. Comparing with the Fourier coefficients of weight two newforms, we found the modular forms $f24B$ and $f120E$.

\subsection{The Galois representation on $H^3_{\et}(X,\QQ_\ell)$}
\label{galX}
The $G_\QQ$-representation on $H^3_\et(X,\QQ_\ell)$ has a summand which is
$H^3_\et(X,\QQ_\ell)^H\cong  H^3_\et(\hY,\QQ_\ell)$, which we discussed in section \ref{comps}, we will assume the conjectured modularity stated there.

Its orthogonal complement $H^3_\et(X,\QQ_\ell)_c$
w.r.t.\ the intersection form was decomposed as a sum of four irreducible $G$-representations of dimension $15,30,45,45$ and each has multiplicity two
(cf. section \ref{decXG}):
$$
H^3_\et(X,\QQ_\ell)_c\,\cong\,V_{15,\ell}\otimes W_{15,\ell}\;\oplus\;
V_{30,\ell}\otimes W_{30,\ell}\;\oplus\;
V_{45,\ell}\otimes W_{45,\ell}\;\oplus\;
V_{45,\ell}'\otimes W_{45,\ell}'.
$$
Each $W_{i,\ell}$ is a $G_{\QQ(i)}$-representation
(but not a $G_\QQ$-representation in general, as the actions of
$G$ and $G_\QQ$ do not commute).
In particular, for $q\equiv 1$ mod 4, the trace of the Frobenius $F_q$ on
$H^3_\et(X,\QQ_\ell)_c$ should be divisible by $15$ in $\ZZ$.
Using a computer we found that
\begin{eqnarray*}
tr(F_q|H^3_\et(X,\QQ_\ell)_c)&=&
tr(F_q|H^3_\et(X,\QQ_\ell))\,-\,tr(F_q|H^3_\et(\hY,\QQ_\ell))\\
&=& \sharp X(\FF_q)-(1+q+q^2+q^3)\,-\,tr(F_q|H^3_\et(\hY,\QQ_\ell))
\end{eqnarray*}
is divisible by $45$ for all such  small $q$.
This leads us to conjecture that $W_{15,\ell}\cong W_{30,\ell}$ as $G_{\QQ(i)}$-representations.
Thus we conjecture:
{\renewcommand{\arraystretch}{1.5}
$$
\begin{array}{rll}
H^3_\et(X,\QQ_\ell)_c &\cong&
(V_{15,\ell}\,\oplus\,V_{30,\ell})\otimes W_{15,\ell}\;\oplus\;
V_{45,\ell}\otimes W_{45,\ell}\;\oplus\;
V_{45,\ell}'\otimes W_{45,\ell}'\\
&\,\stackrel{?}{\cong}_{G_{\QQ(i)}}&\left(W_{15,\ell}\,\oplus\,W_{45,\ell}\,\oplus\,
W_{45,\ell}'\right)^{\oplus 45}.
\end{array}
$$
}
To determine the characteristic polynomial of $F_p$ on
$W_{15,\ell}\oplus W_{45,\ell}\oplus W_{45,\ell}'$
for $p\equiv 1$ mod $4$,
we would have to compute points on $X$ over $\FF_{p^k}$, $k=1,2,3$.
This took too much time,
but we could compute the number of points over
$\FF_p$ and $\FF_{p^2}$ for some small primes.
The characteristic polynomial must be of the form
$$
f_p\,:\quad X^6-s_1X^5+s_2X^4-s_3X^3+ps_2X^2-p^2s_1X+p^3,
$$
with coefficients:
$$
s_1:=t_p,\quad
s_2:=(1/2)t_p^2-(1/2)t_{p^2},\quad
s_3:=(1/6)t_p^3 - (1/2)t_{p^2}t_p+ (1/3)t_{p^3},
$$
where,
$$
t_{p^k}\,:=\,tr(F_p^k|H^3_\et(X,\QQ_\ell)_c(1)),\qquad
|t_{p^k}|\, < \,6\sqrt{p^k}
$$
by Weil's estimate on the eigenvalues of $F_{p^k}$.
For all the primes $1$ mod $4$ from $13$ to $41$ we computed $t_p,t_{p^2}$ and we found a unique integer $n$,
with $|n|< 6\sqrt{p^3}$, such that the polynomial $f_p$ had three quadratic factors upon substituting $t_{p^3}:=n$.
Moreover, comparing the coefficients of the quadratic factors
(of the type $X^2-m_pX +p$ for an integer $m_p$),
we found three elliptic modular forms of weight two whose Fourier coefficients were equal to the coefficients $m_p$. These forms are
$f15C$ (notation as in \cite{modfun4}) and the forms $f24B$ and $f120E$ as in section \ref{comps}. Some Fourier coefficients, denoted by $d_p$, of $f15C$ are given in the table there.

Using this, we next tried to understand the $G_\QQ$-representation on $H^3_\et(X,\QQ_\ell)_c(1))$,
so to find $tr(F_q|H^3_\et)$ also for the $q\equiv 3$ mod 4.
For small $p,k$ we found
$$
tr(F_{p^k}|H^3_\et(X,\QQ_\ell)_c(1))\,=\,\left\{
\begin{array}{rcl} 45(b_p+c_p+d_p)&\mbox{if}&p^k\equiv 1\;\mbox{mod}\,4,\\
9(b_p+c_p+d_p)&\mbox{if}&p^k\equiv 3\;\mbox{mod}\,4.
\end{array}\right.
$$
This leads us to conjecture:
$$
H^3_\et(X,\QQ_\ell)_c(1))\,\stackrel{?}{\cong}_{G_\QQ}\,
W_{c,\ell}^{\oplus 27}\,\oplus\,(\sigma_{1,0,0}\otimes W_{c,\ell})^{\oplus 18},
$$
with the six dimensional $G_\QQ$-representation
$$
W_{c,\ell}\,:=\,V_{f{15C,\ell}}\oplus V_{f{24B,\ell}}\,\oplus\,V_{f{120E,\ell}}
$$
where the $V_{g,\ell}$ are the $\ell$-adic
$G_\QQ$-representations corresponding to the weight two
newforms $f15C$, $f24B$ and $f120E$ respectively.

Putting all conjectures together, one would have the following formula:
$$
\sharp X(\FF_p)\,\stackrel{?}{=}\,1+p+p^2+p^3\,-\,\left\{
\begin{array}{rcl}
(a_p\,+\,p(54b_p\,+\,50c_p\,+\,45d_p)&\mbox{if}&p\equiv 1\;\mbox{mod}\,4,\\
(a_p\,+\,p(18b_p\,+\,23c_p\,+\;9d_p)&\mbox{if}&p\equiv 3\;\mbox{mod}\,4,
\end{array}\right.
$$
we checked that equality holds for all primes $p$ with $7\leq p\leq 500$.

\section{Rational curves and the Abel-Jacobi map for Maschke's Calabi-Yau threefold }\label{lines}

\subsection{Outline}
We show in section \ref{ftl} that Maschke's CY threefold $X$
contains one dimensional families of rational curves.
One such family is denoted by $\cL\rightarrow\widetilde{C}_+$,
it is parametrized by a curve $\widetilde{C}_+$ of genus $33$.
Other families can be obtained by applying the action of $g\in G$ to this family.

From this family one obtains a
morphism of Hodge structures
$$
\phi\,:H^1(\widetilde{C}_+,\ZZ)\,\longrightarrow\,H^3(X,\ZZ).
$$
Similarly one has maps
$\phi_\ell:H^1_\et(\widetilde{C},\QQ_\ell)\rightarrow H^3_\et(X,\QQ_\ell)$ which
are maps of Galois representations (up to Tate twist). To understand $\phi$ and thus also the $\phi_\ell$, we use the Abel-Jacobi map in section \ref{ajm}.

In section \ref{arcp} we study the Galois representation on
$H^1_{\et}(\widetilde{C}_+,\QQ_\ell)$ and relate it, via the map $\phi_{\ell}$,
to the Galois representation on $H^3_{\et}(X,\QQ_\ell)$.

\subsection{Four-tangent lines to Maschke's octic surface}\label{ftl}
To obtain families of rational curves in $X$, we consider curves in $X$ which map isomorphically to lines in $\PP^3$ under the 2:1 map $\pi:X\rightarrow \PP^3$.
A line in $\PP^3$ intersects the branch locus $S$ of $\pi$ in a divisor $D$ of degree $8$, so in general the inverse image of a line is a (hyperelliptic) genus three curve. However, if the line is four tangent to $S$, so $D=2(p_1+\ldots+p_4)$, then the inverse image will
split into two rational curves, each of which maps isomorphically to the line.

It is a classical result that for a general surface of degree $8$ in $\PP^3$ there
are $14752$ lines which are four-tangent to it (\ \cite{Mo}, p.261).
However, $S$ is rather special and in fact it does have positive dimensional families
of four-tangent lines, as we show with some explicit computations.

For $c=(x,y)\in \bA^2$ we consider the line $l_c\subset\bP^3$, with parameter $t$:
$$
l_c\,:\quad (x_0,\,x_1,\,x_2,\,x_3)\,:=\,(x,\,1,\,ty,\,t)
\qquad (c=(x,y)\in \bA^2,\; t\in\CC).
$$
The intersection of $l_c$ with the surface $S$ is defined by the polynomial $f_c$:
$$
l_c\cap S:\quad f_c(t):=F(x,\,1,\,ty\,,t)=0.
$$
The element $g_1\in G$ acts as
$(x_0,x_1,x_2,x_3)\mapsto (x_0,x_1,ix_2,ix_3)$,  with $i^2=-1$, and thus
(as is also easy to verify directly) $F(x_0,x_1,x_2,x_3)=F(x_0,x_1,ix_2,ix_3)$.
Therefore we have $f_c(t)=f_c(it)$ and $f_c(t)$ is actually a polynomial
in $t^4$ of degree $2$.
An explicit computation shows that
$$
f_c(t)\,=\,At^8\,+\,Bt^4\,+\,C
$$
with
$$
A:=y^8 + 14y^4 + 1,\qquad
B:=14(x^4y^4 + x^4 + y^4 + 12x^2y^2 +1),\qquad
C:=x^8 + 14x^4 + 1.
$$

Now we impose that $l_c$ is four-tangent to the surface $S$ by requiring that this polynomial of degree two
in $t^4$ has a double zero. So we restrict the point $(x,y)$ to the algebraic subset of $\bA^2$ defined by
$\Delta=B^2-4AC=0$, where
$$
\Delta\,:=\,(14x^4y^4 + 14x^4 + 168x^2y^2 + 14y^4 + 14)^2\,-\,4( x^8 + 14x^4 + 1)(y^8 + 14y^4 + 1).
$$
For each $c\in\bA^2$ with $\Delta(c)=0$ we have $f_c(t)=(\lambda t^4+\mu)^2$
for some $\lambda,\mu\in\CC$ and thus $l_c$ is a four tangent line.
Therefore $S$ has a positive dimensional family of four-tangent lines.

To understand the locus $\Delta=0$ better, we observe that $\Delta$ is reducible:
$$
\Delta\,=\,48g_+ g_-\qquad(\in \ZZ[x,y])
$$
with polynomials
$$
g_+(x,y)\,:=\,(2y^4 +y^2 +2)x^4-(y^4-24y^2+1)x^2+2y^4 +y^2 +2
$$
and
$$
g_-(x,y)\,:=\,(2y^4 -y^2 +2)x^4+(y^4+24y^2+1)x^2+2y^4 -y^2 +2,
$$
note that $g_-(x,y)=g_+(ix,iy)$.
We define two curves, $C_\pm$, by taking the closure
in $\PP^1\times\PP^1$ of the loci $g_\pm(x,y)=0$ in $\bA^2$.
These curves are isomorphic to the closure of the locus of the lines
parametrized by $g_{\pm}=0$ in the Grassmannian of lines in $\PP^3$.
We checked that $C_\pm$
are smooth curves. As they have bidegree $(4,4)$, their genus is $3^2=9$.
$$
C_+:\quad (2y^4 +y^2v^2 +2v^4)x^4-(y^4-24y^2v^2+v^4)x^2u^2+
2y^4 +y^2v^2 +2v^4\,=\,0\qquad(\subset\,\PP^1_{(x:u)}\times\PP^1_{(y:v)}).
$$

\subsection{Remark: lines in the surface}\label{lins}
Some of the lines $l_c$, parametrized by $c\in C_\pm$, do lie entirely within $S$.
For such a line the coefficients of the degree eight polynomial in $f_c(t)$ must be zero, so:
$$
A\,=\,B\,=\,C\,=\,0.
$$
One verifies that there are $32$ such lines:
$$
l_{(x,y)}\,\subset\, S\quad \Longleftrightarrow\quad
x^8+14x^4+1=0\quad\mbox{and}\quad
y^2=x^2 \;\mbox{or}\;y^2=x^{-2}
$$
and that the corresponding $c=(x,y)$ are the points of intersection of $C_+$ and $C_-$:
$$
l_c\,\subset\, S\qquad \Longleftrightarrow\quad c\, \in \,C_+\cap C_- .
$$
As $C_+$ and $C_-$ are curves of type $(4,4)$ on $\PP^1\times\PP^1$
their intersection number is $4^2+4^2=32$, in particular, their intersection is transversal.
Using the factorization
$$
x^8+14x^4 + 1\,=\,
(x^4 - 2x^3 + 2x^2 + 2x + 1)
(x^4 + 2x^3 + 2x^2 - 2x + 1),
$$
one easily checks that the line $l_3\subset S$ from section \ref{NS(S)} is one of the lines $l_c$ for $c\in C_+\cap C_-$.

\subsection{The double cover $\widetilde{C}_+$ of $C_+$} \label{dccp}
Now we consider the inverse images of the lines $l_c$
with $c\in C_+$ in
the CY threefold $X$. Locally we have
$$
X:\quad w^2\,=\,F(x_0,1,x_2,x_3),\qquad \pi^{-1}(l_c):\quad w^2\,=\,At^8+Bt^4+C,
$$
where $t$ runs over $\CC$.
On $C_+$ we have $B^2=4AC$ and thus, for each $c=(x,y)\in C_+$
with $A(c)\neq 0$ we can rewrite the equation of $ \pi^{-1}(l_c)$ as
$$
4Aw^2\,=\,4A^2t^8+4ABt^4+4AC\,=\,4A^2t^8+4ABt^4+B^2\,=\,
(2At^4+B)^2.
$$
The two irreducible components of the $\pi^{-1}(l_c)$ are thus
defined by the two factors of
$$
(2\sqrt{A}w-(2At^4+B))(2\sqrt{A}w+(2At^4+B)).
$$
Thus the points of the double cover $\widetilde{C}_+$ of $C_+$ defined by $\sqrt{A}$,
where we view $A=y^8 + 14y^4 + 1$ as a rational function on $C_+$,
parametrize rational curves in $X$ which map to lines $l_c$ with $c\in C_+$.
The double cover $\widetilde{C}_+\rightarrow C_+$ is branched over
$4\cdot 8=32$ points, and as ${C}_+$ has genus $9$, $\widetilde{C}_+$ has genus $33$.

\subsection{The Abel-Jacobi map}\label{ajm}
We refer to \cite{V1}, chapter 12, for this section.
We denote by $\cL$ the total space of the family of lines in $X$ parametrized by $\widetilde{C}_+$, it is a surface in the product $\widetilde{C}_+\times X$. In the diagram $f,g$ are the projection maps:
{\renewcommand{\arraystretch}{1.4}
$$
\begin{array}{rcccr}
\widetilde{C}_+\times X&\supseteq&\cL&\stackrel{f}{\longrightarrow}& X\qquad\\
&&\downarrow g&\\
&&\widetilde{C}_+&&
\end{array}\qquad
\cL\,:=\,\{(\widetilde{c},x)\subset \widetilde{C}_+\times X\,:\;x\in\,l_{\widetilde{c}}\,\}.
$$
}
The subvariety $\cL$ has cohomology class, of Hodge type $(2,2)$,
$$
[\cL]\,\in\, H^4(\widetilde{C}_+\times X,\ZZ)\,=\,
\oplus_{i=0}^2 H^i(\widetilde{C}_+,\ZZ)\otimes H^{4-i}(X,\ZZ).
$$
By Poincar\'e duality, $H^1(\widetilde{C}_+,\ZZ)$ is selfdual and thus
$$
H^1(\widetilde{C}_+,\ZZ)\otimes H^{3}(X,\ZZ)\cong Hom_\ZZ(H^1(\widetilde{C}_+,\ZZ),H^3(X,\ZZ)).
$$
As $[\cL]$ has type $(2,2)$, it induces a morphism of
Hodge structures
$$
\phi\,:=\,[\cL]_1\,:H^1(\widetilde{C}_+,\ZZ)\,\longrightarrow\,H^3(X,\ZZ).
$$
This map relates the cohomology of $\widetilde{C}_+$ and $X$.
To understand $\phi$ we use the Abel-Jacobi map.

The (Griffiths) intermediate Jacobian
$$
J(X)\,:=\,( H^{1,2}(X)\oplus H^{0,3}(X))/H^3(X,\ZZ)\;\cong\;
( H^{3,0}(X)^*\oplus H^{2,1}(X)^*)/H^3(X,\ZZ),
$$
where the ${}^*$ indicates the dual complex vector space,
is a $150$-dimensional complex torus.
The Jacobian of $\widetilde{C}_+$ is the abelian variety
$$
J(\widetilde{C}_+)\,:=\,H^{0,1}(\widetilde{C}_+)/H^1(\widetilde{C}_+,\ZZ)\;\cong\;
H^{1,0}(\widetilde{C}_+)^*/H^1(\widetilde{C}_+,\ZZ).
$$
The morphism of Hodge structures $\phi$ corresponds
to the holomorphic map of complex tori
$$
\Phi\,:\;J(\widetilde{C}_+)\,\longrightarrow\,J(X)
$$
which is induced, by the Albanese property of $J(\widetilde{C}_+)$,
by the (holomorphic) Abel-Jacobi map
$$
\Phi_{\widetilde{C}_+}\,:\,\widetilde{C}_+\,\longrightarrow\,J(X),\qquad
\widetilde{c}\longmapsto
\int_\Gamma\;\in\;( H^{3,0}(X)^*\oplus H^{2,1}(X)^*)/H^3(X,\ZZ).
$$
Here one fixes a base point $\widetilde{c}_0\in \widetilde{C}_+$, and $\Gamma\subset X$ is a differentiable 3-chain in $X$ with boundary
$$
\partial \Gamma \,=\,l_{\widetilde{c}}\,-\,l_{\widetilde{c}_0}.
$$
We now point out some easy facts on the domain and the codomain of the Abel-Jacobi map $\Phi_C$ and thus of $\phi$.

The Hodge structure $H^1(\widetilde{C}_+,\QQ)$ has a decomposition into eigenspaces for the covering involution:
$$
H^1(\widetilde{C}_+,\QQ)\,=\,H^1(\widetilde{C}_+,\QQ)_+\,\oplus\,
H^1(\widetilde{C}_+,\QQ)_-,\qquad
H^1(\widetilde{C}_+,\QQ)_+\,\cong\,H^1(C_+,\QQ).
$$
The double cover $\widetilde{C}_+\rightarrow C_+$
induces the projection on the first factor in the eigenspace
decomposition.

A point $c\in C_+$ corresponds to the cycle $\pi^{-1}(l_c)\subset X$.
But as any two lines in $\PP^3$ are fibers of a family of lines
parametrized by a $\PP^1$,
the Abel-Jacobi map restricted to the image of
$J(C_+)\subset J(\widetilde{C}_+)$ is trivial,
since there are no non-constant holomorphic maps from $\PP^1$ to a complex torus.
Hence $\Phi$ factors over the quotient abelian variety $P(\widetilde{C}_+/C_+):=J(\widetilde{C}_+)/J(C_+)$,
the Prym variety of the double cover $\widetilde{C}_+\rightarrow C_+$,
and $\phi$ factors over $H^1(\widetilde{C}_+,\QQ)_-$.

As $\phi$ is a morphism of Hodge structures, we have
$\phi(H^{p,q}(\widetilde{C}_+))\subset H^{p+1,q+1}(X)$.
In section \ref{decXG} we showed that the $G$-action splits the rational Hodge structure on $H^3(X,\QQ)$ into two Hodge substructures:
{\renewcommand{\arraystretch}{1.2}
$$
H^3(X,\QQ)\,=\,H^3(X,\QQ)_t\;\oplus\;H^3(X,\QQ)_a,\qquad
\left\{\begin{array}{rclcrcl}
H^{3,0}(X)_t&\cong&\CC,&\;&H^{2,1}(X)_t&=&0,\\
H^{3,0}(X)_a&\cong&0,&\;&H^{2,1}(X)_a&\cong&\CC^{149},\\
\end{array}\right.
$$
}
\noindent
here $H^3(X,\QQ)_t\cong V_1\otimes W_1$ and $H^3(X,\QQ)_a$
is a direct sum of six summands.
In particular, $H^{1,2}(X)_a=H^{1,2}(X)$.
This implies that  the intermediate Jacobian $J(X)$
has a codimension one subtorus
$$
J_a(X):=H^{1,2}(X)/H^3(X,\ZZ)_a,\qquad
H^3(X,\ZZ)_a\,:=\,H^3(X,\ZZ)\cap H^3(X,\QQ)_a,
$$
which is actually an abelian variety,
polarized by the restriction of the intersection form on $H^3(X,\QQ)$ to $H^3(X,\QQ)_a$.
Moreover, $J_a(X)$ is isogeneous to a product of $149$ elliptic curves (cf.\ section \ref{decXG}).
The image of $\Phi$ is thus contained in $J_a(X)$. Thus the map $\phi$ induces a morphism of Hodge structures:
$$
\phi_a\,:\,H^1(\widetilde{C}_+,\QQ)_-\,\longrightarrow\, H^3(X,\QQ)_a.
$$
In section \ref{ajnt} we show that $\Phi$, and hence $\phi_a$ is non-trivial.

\subsection{Remark}
The surface $S$ contains two $G$-orbits of lines, cf.\ section \ref{NS(S)}.
In Remark \ref{lins} we observed that the lines in the orbit of $l_3$ deform to four tangent lines of $S$. Thus these lines, viewed as cycles with multiplicity two, deform in $X$. On the other hand, we checked that the normal bundle of the line $l_5$ in $X$ is $\cO(-1)^{\oplus 2}$, and thus
these lines are rigid in $X$.
The 1-cycle $l_3-l_5$ might therefore give a non-trivial element in the Griffiths group of $X$.
To check this, one could try to determine the image of this cycle in the intermediate Jacobian $J(X)$,
in particular the component in the isogeny factor
$J_t(X):=H^{0,3}(X)/H^3(X,\ZZ)_t$.
The conjectures of Bloch and Beilinson imply that the rank of the Griffiths group of $X$ over a number field $L$ is the order of vanishing of
the $L$-series of the $G_L$-representation on $W_{1,\ell}$.
As this is (conjectured to be) the restriction of the $G_\QQ$-representation associated to the weight four newform $f120$
(cf.\ section \ref{comps}),
one might be able to compute the order of zero of the L-series in case $L=\QQ(\zeta_{60})$, the field of definition of the 1-cycle.
We refer to \cite{Bloch} for an interesting example involving Schoen's
rigid CY threefold.

In case one has a family of CY-threefolds and a relative family of 1-cycles,
there is an associated inhomogeneous Picard-Fuchs equation.
A nice example, involving the Dwork family of quintic threefolds, appears in \cite{MW}.

\subsection{The Abel-Jacobi map is non-trivial}\label{ajnt}
To show that the Abel-Jacobi map $\Phi$ is non-trivial,
we consider the differential $\rmd \Phi_{\widetilde{C}_+,\widetilde{c}}$
of the map $\Phi_{\widetilde{C}_+}$ which induces $\Phi$
at a point $\widetilde{c}\in\widetilde{C}_+$.
This differential is called the infinitesimal Abel-Jacobi map.
There is a natural map from the deformation space of the pair $(l,X)$, where $l$ is a smooth curve in the Calabi-Yau threefold $X$, to the deformation space of $X$.
If this map is not surjective, then the infinitesimal Abel-Jacobi map is non-zero
(\cite{W}, Lemma 5.1).

The deformations $X_\epsilon$ of $X$ are given by deformations $S_\epsilon$ of Maschke's octic $S$, thus they are defined by $F+\epsilon G$, where $G$ is homogeneous of degree $8$ in $X_0,\ldots,X_3$ and $\epsilon^2=0$.
Considering a rational curve $l_{\widetilde{c}}$
in $X$ which maps onto the line $l_c$,
this curve deforms to $X_\epsilon$ iff the line
$l_c$ deforms to a four tangent line
$l_{c,\epsilon}$ of $S_\epsilon$. As $l_c$ is parametrized by $(x:1:ty:t)$, with $c=(x,y)\in C_+$, any deformation of $l_c$ can be parametrized as
$$
l_{c,\epsilon}:\quad t\,\longmapsto\,(x+\epsilon(a+ct),1,ty+\epsilon(b+dt),t),
$$
for some $a,b,c,d\in\CC$.
Substituting this parametrization in $F+\epsilon G$ one obtains $h(t)+\epsilon k(t)$, for certain polynomials $h,k\in\CC[t]$. As $c\in C_+$, we have $h(t)=f_c(t)^2$ and $l_{c,\epsilon}$ is four tangent to $S_\epsilon$ iff
$f_c(t)^2+\epsilon k(t)\,=\,(f_c(t)+\epsilon m(t))^2$ for some polynomial $m$.
Equivalently, $k(t)=2f_c(t)m(t)$ for some polynomial $m$.
Taking a point $c=(x,y)\in C_+$ with $y=2$ and taking $G=X_0^2X_1^4X_2^2$,
we found that
the coefficient of $t^2$ in $k(t)$ mod $f_c(t)$
is a non-zero constant (i.e.\ independent of $a,b,c,d$), hence there is no polynomial $m(t)$ such that
$k(t)=2f_c(t)m(t)$ and thus the infinitesimal Abel-Jacobi map in ${\widetilde{c}}$
is non-zero.

\section{The geometry and arithmetic of $\widetilde{C}_+$}\label{arcp}

\subsection{Outline}
The map $\phi$ from section \ref{ajm} is defined by the algebraic cycle $[\cL]$ on $\tilde{C}_+\times X$ which is defined over $\QQ$. Hence, after suitable restrictions, we also get maps
$$
\phi_{a,\ell}:H^1_\et(\widetilde{C}_+,\QQ_\ell)_-
\,\longrightarrow\,
H^3_\et(X,\QQ_\ell)_a
$$
which are, up to a Tate twist,
compatible with the action of the Galois group $G_\QQ$.
Thus if $H^1_\et(\widetilde{C}_+,\QQ_\ell)_-$ has a
$G_\QQ$-subrepresentation which does not occur in
$H^3_\et(X,\QQ_\ell)_a$, then that subrepresentation is
mapped to zero by $\phi_{a,\ell}$ and hence also $\phi_a$ has a non-trivial kernel.

Our results indicate that the
$66$ dimensional $G_\QQ$-representation on
$H^1_\et(\widetilde{C}_+,\QQ_\ell)$
might be the direct sum of $33$ two-dimensional $G_\QQ$-representations.
If this is the case, then
the Jacobian of $\widetilde{C}_+$ is isogeneous
to the product of $33$ elliptic curves.

Up to isomorphism, only the Galois representations
associated to six newforms (and their twists by $\sigma_{1,0,0}$)
of weight two
appear in $H^1_\et(\widetilde{C}_+,\QQ_\ell)$.
Up to twist by $\sigma_{1,0,0}$, three of the six $G_\QQ$-representations
in $H^1_\et(\widetilde{C}_+,\QQ_\ell)$ occur only in $H^1_\et(\widetilde{C}_+,\QQ_\ell)_+$ and none of these representations
or their twist appear, even after a Tate twist, in $H^3_\et(X,\QQ_\ell)_a$.
The remaining $G_\QQ$-representations appear only in $H^1(\tilde{C}_+\QQ_\ell)_-$.
All of these do appear, after Tate twist, in $H^3_\et(X,\QQ_\ell)_a$.

In conclusion, we do not  find an obvious obstruction to the surjectivity of $\phi_{a,\ell}$.
We will leave a more detailed study of $\phi_a$ and $\phi_{a,\ell}$ to another occasion.

\subsection{The Galois representation on $H^1_\et(\widetilde{C}_+,\QQ_\ell)$} \label{qtc+}
The genus $33$ curve $\widetilde{C}_+$ is a double cover, defined by
$t^2=A$,  of the genus $9$ curve
$C_+\subset\PP^1\times\PP^1$ defined by
$g_+(x,y)=Px^4-Qx^2+P=0$ where
$$
A\,:=\,y^8 + 14y^4 + 1,\qquad P\,:=\,2y^4 +y^2 +2,\quad Q\,:=\,y^4-24y^2+1.
$$
Hence there is a commutative diagram (see section \ref{dccp})
$$
\begin{array}{rcl} C_+&\longleftarrow&\widetilde{C}_+\\
\downarrow&&\downarrow\\
\PP^1_y&\longleftarrow&C_3,\qquad C_3:\;t^2=A,
\end{array}
$$
where $C_3$ is a hyperelliptic genus three curve.

The covering involution induces the decomposition of $G_\QQ$-representations,
similar to the decomposition of the Hodge structures in section \ref{ajm}:
$$
H^1_\et(\widetilde{C}_+,\QQ_\ell)\,=\,
H^1_\et(\widetilde{C}_+,\QQ_\ell)_+\,\oplus \,H^1_\et(\widetilde{C}_+,\QQ_\ell)_-,
\qquad
H^1_\et(\widetilde{C}_+,\QQ_\ell)_+\,\cong\,H^1_\et({C}_+,\QQ_\ell).
$$
In particular, $H^1_\et({C}_3,\QQ_\ell)\subset H^1_\et(\widetilde{C}_+,\QQ)_-$. As $\sharp C_+(\FF_p)=1-tr(F_p|H^1_\et(\widetilde{C}_+,\QQ_\ell)_+)+p$,
we get:
$$
\sharp\widetilde{C}_+(\FF_p)\,=\,
1-\,tr(F_p|H^1_\et(\widetilde{C}_+,\QQ_\ell))\,+p
\,=\, \sharp {C}_+(\FF_p)\,-\,
tr(F_p|H^1_\et(\widetilde{C}_+,\QQ_\ell)_-).
$$

We computed the cardinality of $\widetilde{C}_+(\FF_p)$ and $C_+(\FF_p)$
for all primes $p$ with $7\leq p\leq 1000$.
The results are consistent with the following formula:
$$
tr(F_p|H^1_\et(\widetilde{C}_+,\QQ_\ell)_-)\,\stackrel{?}{=}\,
(9+3\sigma_{1,0,0}(F_p))b_p+ (5+\sigma_{1,0,0}(F_p))c_p+ (4+2\sigma_{1,0,0}(F_p))d_p
$$
with $b_p,c_p,d_p$ the Fourier coefficients of the newforms $f24B$, $f120E$ and $f15C$,
cf.\ section \ref{comps}.
If this equality holds for all primes $p>5$, or at least for a large set of primes (but we don't know a good bound for this set), then we would have an isomorphism of
$G_\QQ$-representations
$$
H^1_\et(\widetilde{C}_+,\QQ_\ell)_-\,\stackrel{?}{\cong}
{V_{f24B,\ell}}^{\oplus 9}\oplus {V'_{f24B,\ell}}^{\oplus 3}\,\oplus\,
{V_{f120E,\ell}}^{\oplus 5}\oplus {V'_{f120E,\ell}}\,\oplus\,
{V_{f15C,\ell}}^{\oplus 4}\oplus {V'_{f15C,\ell}}^{\oplus 2},
$$
where the twist by $\sigma_{1,0,0}$ of a $G_\QQ$-representation is denoted by
$$
V'_{\ast}\,:=\,V_{\ast}\otimes\sigma_{1,0,0}.
$$
Note that the three representations and their twists in this decomposition
appear, after a Tate twist, in $H^3(X,\QQ_\ell)_a$.
See Remark \ref{remc} for a possible geometric approach to the conjectured decomposition, with the observations made there one can actually prove that
$$
H^1_\et({C}_3,\QQ_\ell)\,\cong \,
{V_{f24B,\ell}}^{\oplus 2}\oplus {V'_{f24B,\ell}}.
$$

Similarly, we conjecture (and we checked equality of the traces of $F_p$
for primes $p$ with $5<p<1000$):
$$
H^1_\et({C}_+,\QQ_\ell)\,\stackrel{?}{\cong}
{V_{f210,\ell}}^{\oplus 3}\,\oplus\,
{V_{f840,\ell}}^{\oplus 2}\oplus {V'_{f840,\ell}}\,\oplus\,
V_{f1680,\ell}\oplus {V'_{f1680,\ell}}^{\oplus 2} ,
$$
where the $f_N$ are newforms of weight $2$ on $\Gamma_0(N)$ characterized
by the Fourier coefficients in the table below.
{\renewcommand{\arraystretch}{1.2}
$$
\begin{array}{|r|r|r|r|r|r|r|r|r|r|r|r|r|r|r|}\hline
&p &11&13&17&19&23&29&31&37&\ldots&79&83&89&97 \\
\hline
f_{210},k=2 &b'_p&4 &-2& 2&-4 &-8&6&-8&-2&\ldots&0&12&2&10\\
\hline
f_{840},k=2&c'_p&-4&-2&2&-4&0&-10&0&6&\ldots&8&-4&10&10\\
\hline
f_{1680},k=2&d'_p&-4&-2&2&-4&0&-2&-8&-2&\ldots&8&-4&2&-14\\
\hline
\end{array}
$$
}

\

Note that none of these Galois  representations occurs
(even after Tate twist) in $H^3_\et(X,\QQ_\ell)$.

\subsection{Remark: quotients of $\widetilde{C}_+$}\label{remc}
The curve $\tilde{C}_+$ has a big automorphism group, which one can use to decompose the cohomology and to find quotient curves, of lower genus, of $\tilde{C}_+$.
One might hope these automorphisms, and further automorphisms of the quotient curves not induced automorphisms of $\tilde{C}_+$, might suffice to prove the conjectural decomposition of the \'etale cohomology. However, we did not succeed in carrying this out, but the partial results we obtained were quite helpful in finding the conjectured decomposition.

The group $(\ZZ/2\ZZ)^2$, generated by the involutions
$$
\iota_1:\;(x,y,t)\,\longmapsto\,(-x,y,t),\qquad
\iota_2:\;(x,y,t)\longmapsto\, (x,y,-t).
$$
 on $\widetilde{C}_+$. From this one finds the following quotient curves.

The curve $\overline{C}_+=\widetilde{C}_+/(\ZZ/2\ZZ)^2$ is a double cover
$\overline{C}_+\rightarrow \PP^1_y$ defined by $Px^2-Qx+P=0$.
In particular, it is a hyperelliptic curve of genus $3$, defined by $s^2=Q^2-4P^2$.

The quotient of $\widetilde{C}_+$ by the product $\iota_1\iota_2$
of these involutions
is a curve $C_{17}$ of genus $17$ which is a double cover of $\overline{C}_3$.

The curve $C_{13}:=\widetilde{C}_+/\iota_2$ has genus $13$, it is a double cover of both
$\overline{C}_+$ and of $C_3$, each of these covers is branched in $2\cdot 8=16$ points.
This curve is again a $(\ZZ/2\ZZ)^2$-cover $\PP^1_y$, two quotients by involutions are
$\overline{C}_+$ and $C_3$, the third quotient by an involution is
a 2:1 cover of $\PP^1_y$ branched over the $8+8$ ramification points of the other two double covers, so it is a curve $C_7$ of genus 7 with equation
$u^2=A(Q^2-4P^2)=
(y^8+14y^4+1)
(y^4-24y^2+1)^2-4(2y^4+y^2+2))
$.
Thus $C_7$ has an involution
$(u,y)\mapsto (u,-y)$, with quotient curve of genus three.
Another involution on $C_7$, which is fixed point free,
is $(u,y)\mapsto (-u,-y)$ and it has a genus four quotient. Etcetera.
Using similar involutions, one can find genus one quotients of $C_3$ which
lead to the decomposition of $H^1_\et(C_3,\QQ_\ell)$ given in the previous section.

\section{Trace formulas}\label{trfo}
\subsection{} In this section we prove the formulas in Proposition \ref{propch}.
The first two, which give the Euler characteristic and the dimension $h^n(X)_{pr}$
of the primitive cohomology group of
a hypersurface in projective space, are well known. The second is a formula of Ch\^{e}nevert \cite{Ch} which determines the trace of an automorphism of a hypersurface on the primitive cohomology. We derive it directly from the Lefschetz trace formula.
The last is an easy generalization of Ch\^{e}nevert's formula to cyclic ramified covers of projective space,
like Maschke's CY. These trace formulas are remarkable and easy to use, since they involve only the degree of $X$ and  the dimensions of certain eigenspaces of $\sigma$, but not the specific geometry of $X\subset\PP^{n+1}$.

\subsection{Proposition}\label{propch}
Let $X$ be a smooth hypersurface of degree $d$ and dimension $n$ in $\PP^{n+1}$
defined by an equation $F=0$. Let $\sigma:\CC^{n+2}\rightarrow\CC^{n+2}$ be a linear map such that $F(\sigma(x))=F(x)$. Let $r$ be a divisor of $d$ and let $Y\rightarrow \PP^{n+1}$
be the cyclic $r$:1 cover of $\PP^{n+1}$ branched along $X$ and let $\widetilde{\sigma}\in Aut(Y)$ be the automorphism induced by $\sigma$.
Then we have:
\begin{enumerate}
\item
The Euler characteristic of  $X$ is given by:
$$
\chi(X)\,=\,n+2\,+\frac{1}{d}((1-d)^{n+2}-1),\qquad
h^n_{pr}(X)\,:=\,\dim H^n(X,\QQ)_{pr}\,=\,
\frac{(-1)^n}{d}((1-d)^{n+2}+d-1).
$$
\item (Ch\^{e}nevert's formula \cite{Ch})
For $\alpha\in\CC$ such that $\alpha^d=1$,
let $m_\alpha$ be the multiplicity of the eigenvalue $\alpha$ of $\sigma$, where we put $m_\alpha=0$ if $\alpha$ is not an eigenvalue of $\sigma$. Then the trace of automorphism of $X$ induced by $\sigma$
on the primitive cohomology of $X$ is the following sum over
all the $d$-th roots of unity:
$$
tr(\sigma^*|H^n(X,\QQ)_{pr})\,=\,\frac{(-1)^n}{d}
\sum_{\alpha^d=1}(1-d)^{m_\alpha}.
$$
\item
With the notations as above, we have:
$$
tr(\widetilde{\sigma}^*|H^{n+1}(Y,\QQ)_{pr})\,=\,\frac{(-1)^{n+1}}{d}\left(
\sum_{\alpha^d=1}
\mbox{$$}(1-d)^{m_{\alpha}}
\;-\,
r\sum_{\gamma^{d/r}= 1}
\mbox{$$}(1-d)^{m_{\gamma}}\right).
$$
\end{enumerate}

\ts \begin{enumerate}\item
The first formula can be obtained from the Gauss-Bonnet formula $\chi(X)=c_n(T_X)$. The normal bundle sequence implies that
$$
c(T_X)=\frac{c(T_{{\PP^{n+1}|X})}}{c(\cO(d))}\,=\,
\frac{(1+h)^{n+2}}{1+dh}\,=\,
\left( \sum_{i=0}^{n+2} \binom{n+2}{i}h^i\right)
\left(\sum_j (-d)^jh^j\right)
$$
where $h$ is the hyperplane class on $X$, so $h^n=d$ and $h^i=0$ for $i>n$.
In degree $n$ we get:
$$
c_n(T_X)\,=\,d\sum_{k=0}^n\binom{n+2}{k}(-d)^{n-k}.
$$
Comparing this with
$$
(1-d)^{n+2}\,=\,d^2\left(\sum_{k=0}^n \binom{n+2}{k}(-d)^{n-k}\right)
\,-\,d(n+2)+1
$$
gives the formula for $\chi(X)$.

By Lefschetz's hyperplane theorem, one has
$H^{2i}(X,\QQ)\cong H^{2i}(\PP^{n+1},\QQ)\cong \QQ$
if $2i\neq n$ and $H^{2i+1}(X,\QQ)=0$ if $2i+1\neq n$.
Thus the non-primitive cohomology of $X$ contributes $d+1$ to the Euler characteristic, the primitive cohomology is concentrated in $H^{n}(X,\QQ)$, and its dimension is $h^n_{pr}(X)=\chi(X)-(n+1)$. Hence
$$
h^n_{pr}(X)\,=\,(-1)^n(1+\frac{1}{d}((1-d)^{n+2}-1))\,=\,
\frac{(-1)^n}{d}((1-d)^{n+2}+d-1).
$$

\item
An eigenvalue of $\sigma$ is denoted by
$\beta$ and its multiplicity by $m_\beta$, so that $\sum m_\beta=n+2$.
Let $\PP_\beta\subset \PP^{n+1}$ be the projectivization of the eigenspace of $\sigma$ with eigenvalue $\beta$, it has dimension $m_\beta-1$.
Let  $X_\beta:=X\cap\PP_\beta$, then
$X_\beta=\PP_\beta\subset X$ if $\beta^d\neq 1$ and else $X_\beta$
is a smooth subvariety of $\PP_\beta$ of dimension $m_\beta-2$ and degree $d$
(cf.\ \cite{Ch}, Lemma 2.3).

The Lefschetz fixed point formula for an automorphism $\sigma$ of $X$ is:
$$
\sum(-1)^itr(\sigma^*|H^i(X,\QQ))\,=\,\sum_j \chi(X_j),
$$
where the right hand sum is over the components of the fixed point set of
$\sigma$ in $X$ (c.f.\ \cite{Peters}).
As we observed above, the right hand side can be written as a sum over the eigenvalues $\beta$ of $\sigma$:
$$
\sum_j \chi(X_j)\,=\,\sum_{\beta}\chi(X_\beta)\,=\,
\sum_{\beta^d\neq 1} \chi(\PP_\beta)\,+\,\sum_{\beta^d=1}\chi(X_\beta)
$$
where the $X_\beta$ are now all smooth hypersurfaces of degree $d$ in $\PP_\beta\cong\PP^{m_\beta-1}$. Thus we get:
$$
\sum_j \chi(X_j)\,=\,\sum_{\beta^d\neq 1} m_\beta\;+
\,\sum_{\beta^d=1}m_\beta +\mbox{$\frac{1}{d}$}((1-d)^{m_\beta}-1)
$$
which, using $\sum m_\beta=n+2$, simplifies to
$$
\sum_j \chi(X_j)\,=\,n+2\,+\,
\mbox{$\frac{1}{d}$}\sum_{\beta^d=1}(1-d)^{m_\beta}-1).
$$
If $\alpha^d=1$ and $\alpha$ is not an eigenvalue of $\sigma$, then $m_\alpha=0$ and thus $(1-d)^{m_\alpha}-1=0$. There are $d$ complex numbers with $\alpha^d=1$ and so:
$$
\sum_j \chi(X_j)\,=\,n+2\,+\,
\sum_{\alpha^d=1}\mbox{$\frac{1}{d}$}((1-d)^{m_\alpha}-1)\,=\,
n+1\,+\,\mbox{$\frac{1}{d}$}\sum_{\alpha^d=1}(1-d)^{m_\alpha}.
$$
As $\sigma^*$ is trivial on the non-primitive cohomology, the left hand side is
$$
\sum(-1)^itr(\sigma^*|H^i(X,\QQ))\,=\,n+1\,+\,(-1)^ntr(\sigma^*|H^n((X,\QQ)_{pr}),
$$
hence Ch\^{e}nevert's formula follows.

\item
Let $r\in\ZZ_{\geq 0}$ divide $d$ and let $Y$ be the  $r$:1 cover of $\PP^{n+1}$ branched along $X$. The variety $Y$ has a natural embedding in the total space of the line bundle
$\cO(d/r)$ over $\PP^{n+1}$, where it is defined by $t^r=F$.
The automorphism $\sigma$ of $\PP^{n+1}$ lifts to an automorphism $\widetilde{\sigma}$
of the tautological line bundle $\cO(-1)\subset\PP^{n+1}\times\CC^{n+2}$,
its action is induced by $(\sigma,\sigma)$ on $\CC^{n+2}\times\CC^{n+2}$.
This induces an action, again denoted by $\widetilde{\sigma}$,
on $\cO(d/r)$ such that if $x\in\CC^{n+2}$ is an eigenvector of $\sigma$ with eigenvalue $\beta$, then $\widetilde{\sigma}$ acts as scalar multiplication by $\beta^{-d/r}$ on the fiber of $ \cO(d/r)$ over $<x>\in\PP^{n+1}$.

Let $y\in Y$ be a fixed point of $\widetilde{\sigma}$, then its image $x\in\PP^{n+1}$ is a fixed point for $\sigma$, hence $x\in \PP_\beta$ for some eigenvalue $\beta$ of $\sigma$.
In case $\beta^d\neq 1$, $\PP_\beta\subset X$, hence $F(x)=0$
and $x=y$ in the total space of $\cO(d/r)$.
Conversely, any point $x\in \PP_\beta\subset Y$ is a fixed point of $\widetilde{\sigma}$.
In case $\beta^d= 1$, there are two possibilities: if $x\in X_\beta=X\cap\PP_\beta$,
then $x=y$ and conversely any point in $X_\beta\subset Y$
is a fixed point of $\widetilde{\sigma}$.
If $x\notin X_\beta$, then we must have
$\beta^{-d/r}=1$, in that case all points in the preimage $Y_\beta$ of $\PP_\beta$ are fixed points of $\widetilde{\sigma}$.
The variety $Y_\beta$ is a smooth $r$:1 cover of $\PP_\beta$ branched along $X_\beta$.
Thus the fixed point set of $\widetilde{\sigma}$ in $Y$ is:
$$
Y^{\widetilde{\sigma}}\,=\,\left(\coprod_{\beta^{d/r}\neq 1} X_\beta\right)
\;\coprod\;
\left(\coprod_{\beta^{d/r}= 1} Y_\beta\right).
$$

The Euler characteristic of the fixed point set $Y^{\widetilde{\sigma}}$ is thus given by
$$
\chi(Y^{\widetilde{\sigma}})\,=\,\sum_{\beta^d\neq 1} \chi(\PP_\beta)\,+\,
\sum_{\beta^d=1, \beta^{d/r}\neq 1} \chi(X_\beta)\,+\,
\sum_{\beta^{d/r}= 1} \chi(Y_\beta),
$$
where the $X_\beta$ which appear are smooth hypersurfaces of degree $d$ in $\PP_\beta\cong\PP^{m_\beta-1}$. The Hurwitz formula for the branched $r$:1 cover $Y_\beta\rightarrow \PP_\beta$ gives:
$$
\chi(Y_\beta)\,=\,r\chi(\PP_\beta)\,-\,(r-1)\chi(X_\beta)\,=\,
rm_\beta\,-\,(r-1)(m_\beta\,+\,\mbox{$\frac{1}{d}$}((1-d)^{m_\beta}-1)).
$$
Therefore we get:
$$
\chi(Y^{\widetilde{\sigma}})\,=\,\sum_{\beta^d\neq 1} m_{\beta}\;+\,
\sum_{\beta^d=1, \beta^{d/r}\neq 1}
m_{\beta}+\mbox{$\frac{1}{d}$}((1-d)^{m_{\beta}}-1)
\;+\,
\sum_{\beta^{d/r}= 1} m_{\beta}-
\mbox{$\frac{r-1}{d}$}((1-d)^{m_{\beta}}-1),
$$
where we sum over all eigenvalues $\beta$ of $\sigma$.
As $\sum m_{\beta}=n+2$ we obtain:
$$
\chi(Y^{\widetilde{\sigma}})\,=\,n+2\,+\,\mbox{$\frac{1}{d}$}
\sum_{\beta^d=1, \beta^{d/r}\neq 1}
((1-d)^{m_{\beta}}-1)
\;-\,
\mbox{$\frac{r-1}{d}$}\sum_{\beta^{d/r}= 1}
((1-d)^{m_{\beta}}-1).
$$
Since $(1-d)^{m_{\alpha}}-1=0$ if $m_\alpha=0$, we can rewrite the formula with a sum over all $d$-th roots of unity $\alpha$ and a sum over all $d/r$-th roots of unity $\gamma$:
$$
\chi(Y^{\widetilde{\sigma}})\,=\,n+2\,+\,
\frac{1}{d}\sum_{\alpha^d=1}
\mbox{$$}((1-d)^{m_{\alpha}}-1)
\;-\,
\frac{r}{d}\sum_{\gamma^{d/r}= 1}
\mbox{$$}((1-d)^{m_{\gamma}}-1),
$$
which simplifies to
$$
\chi(Y^{\widetilde{\sigma}})\,=\,n+2\,+\,
\frac{1}{d}\sum_{\alpha^d=1}
\mbox{$$}(1-d)^{m_{\alpha}}
\;-\,
\frac{r}{d}\sum_{\gamma^{d/r}= 1}
\mbox{$$}(1-d)^{m_{\gamma}}.
$$

As any automorphism $\widetilde{\sigma}$ is the identity on the non-primitive cohomology of $Y$, we find that
$$
tr(\widetilde{\sigma}|H^{n+1}(Y)_{pr})\,=\,\frac{(-1)^{n+1}}{d}\left(
\sum_{\alpha^d=1}
\mbox{$$}(1-d)^{m_{\alpha}}
\;-\,
r\sum_{\gamma^{d/r}= 1}
\mbox{$$}(1-d)^{m_{\gamma}}\right).
$$
This concludes the proof of the proposition.
\end{enumerate}
\qed

\subsection{Remark}
In case $r=d$, $Y$ is a hypersurface in the total space of $\cO(1)$, which
can be identified with $\PP^{n+2}-\{(0:\ldots:0:1)\}$. The bundle projection to $\PP^{n+1}$ is given by $(x_0:\ldots:x_{n+1}:x_{n+2})\mapsto (x_0:\ldots:x_{n+1})$ and $Y$ is defined by $F(x_0,\ldots,x_{n+1})=x_{n+2}^d$.
The action of $\widetilde{\sigma}$ is as $\sigma$ on $x_0,\ldots,x_{n+1}$ and is trivial on $x_{n+2}$. Hence the eigenspaces $\PP_\beta$ of $\sigma$ and $\widetilde{\sigma}$ are the same if $\beta\neq 1$ and the dimension of $\PP_1$ increases by one. One easily verifies that our generalization of Ch\^{e}nevert's formula in this case gives the original formula.

\end{document}